\documentclass[final]{article}

\usepackage{amsmath}
\usepackage{amsfonts}
\usepackage{color}
\usepackage{graphicx}
\usepackage{subfig}
\usepackage{caption}
\usepackage[margin=1in,headsep=.2in]{geometry}
\usepackage{hyperref}
\hypersetup{
    colorlinks=true,
    linkcolor=blue,
    citecolor=blue,
    filecolor=magenta,
    urlcolor=cyan,
}

\usepackage{algorithmic,algorithm}
\usepackage{enumitem}
\usepackage{comment}
\usepackage{enumitem}
\def\zerobold{\boldsymbol{0}}

\newlist{Assumption}{enumerate}{1}
\setlist[Assumption]{label=A\arabic*}

\def\proof{\par{\it Proof}. \ignorespaces}

\def\endproof{\vbox{\hrule height0.6pt\hbox{%
   \vrule height1.3ex width0.6pt\hskip0.8ex
   \vrule width0.6pt}\hrule height0.6pt
  }}

\newtheorem{theorem}{Theorem}

\newcommand{\bmat}[1]{\begin{bmatrix}#1\end{bmatrix}} 
\definecolor{Blue}{rgb}{0,0,1}
\definecolor{Red}{rgb}{1,0,0}
\definecolor{Green}{rgb}{0,1,0}
\newcommand{\overbar}[1]{\mkern 1.5mu\overline{\mkern-1.5mu#1\mkern-1.5mu}\mkern 1.5mu}

\newcommand{\paramSymbol}{\mu}
\newcommand{\param}{\boldsymbol{\paramSymbol}}
\newcommand{\paramk}[1]{\param_{#1}}
\newcommand{\nparam}{n_{\paramSymbol}}
\newcommand{\statemat}{\boldsymbol{A}}
\newcommand{\reducedstatemat}{\hat{\statemat}}
\newcommand{\ststatemat}{\statemat^{\text{st}}}
\newcommand{\reducedststatemat}{\reducedstatemat^{\text{st}}}
\newcommand{\reducedststatematk}[2]{\reducedststatemat_{(#1,#2)}}
\newcommand{\basis}{\boldsymbol{\Phi}}
\newcommand{\sconstantserror}{M}
\newcommand{\sconstantserrork}[1]{\sconstantserror_{#1}}
\newcommand{\stabilityconstant}{L}
\newcommand{\stabilityconstantk}[1]{\stabilityconstant_{#1}}
\newcommand{\projmat}{\boldsymbol{P}}
\newcommand{\projerrormat}{\boldsymbol{D}}
\newcommand{\timeintegratorop}{\boldsymbol{T}}
\newcommand{\timeintegratoropk}[1]{\timeintegratorop_{#1}}
\newcommand{\remainvec}{\boldsymbol{q}}
\newcommand{\remainveck}[1]{\remainvec_{#1}}
\newcommand{\serrorboundcoef}{\beta}
\newcommand{\serrorboundcoefk}[1]{\serrorboundcoef_{#1}}
\newcommand{\soperatornorm}{\alpha}
\newcommand{\soperatornormk}[1]{\soperatornorm_{#1}}
\newcommand{\spatialbasis}{\basis_{\text{s}}}

\newcommand{\spacetimebasis}{\basis_{\text{st}}}
\newcommand{\basisvec}{\boldsymbol{\phi}}
\newcommand{\spacetimebasisvec}{\basisvec^{\text{st}}}
\newcommand{\spacetimebasisveck}[1]{\spacetimebasisvec_{#1}}
\newcommand{\spatialbasisvec}{\basisvec^{\text{s}}}
\newcommand{\spatialbasisveck}[1]{\spatialbasisvec_{#1}}
\newcommand{\temporalbasis}{\basis_{\text{t}}}
\newcommand{\temporalbasisi}[1]{\temporalbasis^{#1}}
\newcommand{\temporalbasisvec}{\basisvec^{\text{t}}}
\newcommand{\temporalbasisveck}[2]{\temporalbasisvec_{#1#2}}
\newcommand{\basiselement}{\phi}
\newcommand{\temporalbasiselement}{\basiselement^{\text{t}}}
\newcommand{\temporalbasiselementk}[3]{\temporalbasiselement_{#1#2,#3}}
\newcommand{\diagonaltemporalbasis}[2]{\boldsymbol{D}_{#1}^{#2}}
\newcommand{\identity}{\boldsymbol{I}}
\newcommand{\identityk}[1]{\identity_{#1}}
\newcommand{\zero}{\boldsymbol{0}}
\newcommand{\inputmat}{\boldsymbol{B}}
\newcommand{\reducedinputmat}{\hat{\inputmat}}
\newcommand{\outputmat}{\boldsymbol{C}}
\newcommand{\reducedoutputmat}{\hat{\outputmat}}
\newcommand{\inputvec}{\boldsymbol{f}}
\newcommand{\stinputvec}{\inputvec^{\text{st}}}
\newcommand{\reducedstinputvec}{\hat{\inputvec}^{\text{st}}}
\newcommand{\reducedstinputveck}[1]{\reducedstinputvec_{(#1)}}
\newcommand{\outputvec}{\boldsymbol{y}}
\newcommand{\outputveck}[1]{\outputvec^{(#1)}}
\newcommand{\resSymbol}{r}
\newcommand{\res}{\boldsymbol{\resSymbol}}
\newcommand{\resk}[1]{\res^{(#1)}}
\newcommand{\stres}{\res^{\text{st}}}
\newcommand{\approxres}{\tilde{\res}}
\newcommand{\approxresk}[1]{\approxres^{(#1)}}
\newcommand{\stateSymbol}{u}
\newcommand{\stateSnapshot}{\boldsymbol{U}}
\newcommand{\temporalSnapshot}{\boldsymbol{\Upsilon}}
\newcommand{\temporalSnapshotk}[1]{\temporalSnapshot_{#1}}
\newcommand{\stateSnapshotk}[1]{\stateSnapshot_{#1}}
\newcommand{\state}{\boldsymbol{\stateSymbol}}
\newcommand{\ststate}{\state^{\text{st}}}
\newcommand{\statek}[1]{\state^{(#1)}}
\newcommand{\diffstatek}[1]{\delta \state^{(#1)}}
\newcommand{\diffststate}{\delta \ststate}
\newcommand{\initialstate}{\state_0}
\newcommand{\reducedinitialstate}{\hat{\state}_0}
\newcommand{\stinitial}{\initialstate^{\text{st}}}
\newcommand{\reducedstinitial}{\hat{\state}_0^{\text{st}}}
\newcommand{\reducedstinitialk}[1]{\hat{\state}_{0,#1}^{\text{st}}}
\newcommand{\approxstate}{\tilde{\state}}
\newcommand{\approxstatek}[1]{\approxstate^{(#1)}}
\newcommand{\approxststate}{\approxstate^{\text{st}}}
\newcommand{\refstate}{\state_{\text{ref}}}
\newcommand{\reducedrefstate}{\hat{\state}_{\text{ref}}}
\newcommand{\reducedstate}{\hat{\state}}
\newcommand{\reducedststate}{\reducedstate^{\text{st}}}
\newcommand{\reducedstatek}[1]{\reducedstate^{(#1)}}
\newcommand{\timevar}{t}
\newcommand{\timevark}[1]{\timevar_{(#1)}}
\newcommand{\initialtime}{0}
\newcommand{\finaltime}{T}
\newcommand{\timestep}{\Delta \timevar}
\newcommand{\timestepk}[1]{\timestep^{(#1)}}
\newcommand{\inputveck}[1]{\inputvec^{(#1)}}

\newcommand{\position}{\boldsymbol{r}}
\newcommand{\energy}{E}
\newcommand{\energyin}{\energy'}
\newcommand{\speed}{\nu}
\newcommand{\crosssection}{\sigma}
\newcommand{\scatter}{\crosssection_{\text{s}}}
\newcommand{\direction}{\boldsymbol{\Omega}}
\newcommand{\directionin}{\boldsymbol{\Omega}'}
\newcommand{\source}{q}

\newcommand{\leftsingularmat}{\boldsymbol{W}}
\newcommand{\rightsingularmat}{\boldsymbol{V}}
\newcommand{\leftsingularmatk}[1]{\leftsingularmat_{#1}}
\newcommand{\leftsingularmattemporal}{\boldsymbol{\Lambda}}
\newcommand{\leftsingularmattemporalk}[1]{\leftsingularmattemporal_{#1}}
\newcommand{\rightsingularmattemporal}{\boldsymbol{\Psi}}
\newcommand{\rightsingularmattemporalk}[1]{\rightsingularmattemporal_{#1}}
\newcommand{\rightsingularmatk}[1]{\rightsingularmat_{#1}}
\newcommand{\leftsingularvec}{\boldsymbol{w}}
\newcommand{\rightsingularvec}{\boldsymbol{v}}
\newcommand{\leftsingularveck}[1]{\leftsingularvec_{#1}}
\newcommand{\rightsingularveck}[1]{\rightsingularvec_{#1}}
\newcommand{\singularvalmat}{\boldsymbol{\Sigma}}
\newcommand{\singularvalmatk}[1]{\singularvalmat_{#1}}
\newcommand{\middlesingularvalmat}{\overbar{\singularvalmat}}
\newcommand{\basisrank}{r}
\newcommand{\basisrankk}[1]{\basisrank_{#1}}
\newcommand{\maxrank}{r_{\text{max}}}

\newcommand{\middlerightsingularmat}{\overbar{\rightsingularmat}}
\newcommand{\updatedBasisVector}{\boldsymbol{j}}
\newcommand{\projerr}{p}
\newcommand{\reducedsnapshotvec}{\boldsymbol{\ell}}
\newcommand{\middlemat}{\boldsymbol{Q}}

\newcommand{\middleleftsingularmat}{\overbar{\leftsingularmat}}
\newcommand{\tol}{\epsilon}
\newcommand{\svdtol}{\tol_\text{SVD}}
\newcommand{\svtol}{\tol_\text{SV}}
\newcommand{\singularvaluesSymbol}{\sigma}
\newcommand{\singularvalues}{\boldsymbol{\singularvaluesSymbol}}
\newcommand{\singularvaluesk}[1]{\singularvalues_{#1}}
\newcommand{\assign}{\leftarrow}
\newcommand{\singularvaluesArg}[1]{\singularvaluesSymbol_{#1}}
\newcommand{\projected}{\boldsymbol \ell}
\newcommand{\diag}[1]{\text{diag}(#1)}
\newcommand{\Qmat}{\overbar{\boldsymbol{Q}}}
\newcommand{\Rmat}{\boldsymbol R}
\newcommand{\qr}[1]{ QR\{#1\} }

\newcommand{\snapshotvec}{\state}

\newcommand{\NN}{\mathbb{N}}
\newcommand{\RR}[1]{\ensuremath{\mathbb{R}^{ #1 }}}

\newcommand{\dummyvec}{\boldsymbol{a}}

\newcommand{\paramspace}{\Omega_{\paramSymbol}}
\newcommand{\fontDiscrete}{\mathcal}
\newcommand{\paramsample}{{\fontDiscrete P}}

\newcommand{\nspace}{N_s}
\newcommand{\nreducedspace}{n_s}
\newcommand{\nreducedtime}{n_t}
\newcommand{\ntime}{N_t}
\newcommand{\ninput}{N_i}
\newcommand{\noutput}{N_o}

\newcommand{\range}[1]{\text{range}(#1)}

\newcommand{\natNo}{\NN} 
\newcommand{\nat}[1]{\natNo(#1)}
\newcommand{\rddots}{\cdot^{\cdot^{\cdot}}}

\newcommand{\tensor}{\otimes}

\newcommand{\SCB}{\textsc{SCB}}

\newcommand{\reducedtimeindex}{j}
\newcommand{\reducedspaceindex}{i}
\providecommand{\keywords}[1]
{
  \small
  \textbf{\textit{Keywords---}} #1
}

\title{Space--time reduced order model for large-scale linear dynamical
systems with application to Boltzmann transport problems}

\author{Youngsoo Choi\footnote{Correspondence to: Lawrence Livermore National
 331 Laboratory, 7000 East Ave, Livermore, CA 94550, USA. E-mail:
 choi15@llnl.gov}, Peter Brown, Bill Arrighi, Robert Anderson
 \vspace{6pt}
\\Lawrence Livermore National Laboratory\footnote{Lawrence 
Livermore National Laboratory is operated by Lawrence Livermore
National Security, LLC, for the U.S. Department of Energy, National Nuclear
Security Administration under Contract DE-AC52-07NA27344 and LLNL-JRNL-791966.}
}
\date{}

\begin{document}

\maketitle

\begin{abstract}
A classical reduced order model for dynamical problems involves spatial
  reduction of the problem size. However, temporal reduction accompanied by the
  spatial reduction can further reduce the problem size without losing accuracy
  much, which results in a considerably more speed-up than the spatial reduction
  only. Recently, a novel space--time reduced order model for dynamical problems
  has been developed \cite{choi2019space}, where the space--time reduced order
  model shows an order of a hundred speed-up with a relative error of $10^{-4}$
  for small academic problems. However, in order for the method to be applicable
  to a large-scale problem, an efficient space--time reduced basis construction
  algorithm needs to be developed. We present incremental space--time reduced
  basis construction algorithm.  The incremental algorithm is fully parallel and
  scalable. Additionally, the block structure in the space--time reduced basis
  is exploited, which enables the avoidance of constructing the reduced
  space--time basis.  These novel techniques are applied to a large-scale
  particle transport simulation with million and billion degrees of freedom.
  The numerical example shows that the algorithm is scalable and practical.
  Also, it achieves a tremendous speed-up, maintaining a good accuracy. Finally,
  error bounds for space-only and space--time reduced order models are derived.
\end{abstract}

\keywords{Space--time reduced order model, incremental singular value
decomposition, Boltzmann transport equations, linear dynamical systems, block
structure, optimal projection, proper orthogonal decomposition}

\section{Introduction}\label{sec:intro}

    Many computational models for physics simulations are formulated as linear
    dynamical systems. Examples of linear dynamical systems include the
    computational model for the signal propagation and interference in electric
    circuits, storm surge prediction models before an advancing hurricane,
    vibration analysis in large structures, thermal analysis in various media,
    neuro-transmission models in the nervous system, various computational
    models for micro-electro-mechanical systems, and various particle transport
    simulations.  Depending on the complexity of geometries and desirable
    fidelity level, these problems become easily large-scale problems. For
    example, the Boltzmann Transport Equation (BTE) has seven independent
    variables, i.e., three spatial variables, two directional variables, one
    energy variable, and one time variable. It is not hard to see that the BTE
    can easily lead to a high dimensional discretized problem.  Additionally,
    the complex geometry (e.g., reactors with thousands of pins and shields) can
    lead to a large-scale problem.  As an example, a problem with 20 angular
    directions, a cubit spatial domain of 100 x 100 x 100 elements, 16 energy
    groups, and 100 time steps leads to 32 billion unknowns.  The large-scale
    hinders a fast forward solve and prevents the multi-query setting problems,
    such as uncertainty quantification, design optimization, and parameter
    study, from being tractable.  Therefore, developing a Reduced Order Model
    (ROM) that accelerates the solution process without losing much accuracy is
    essential. 

    There are several model order reduction approaches available for linear
    dynamical systems: (i) Balanced truncation \cite{mullis1976synthesis,
    moore1981principal} in control theory community is the most famous one. It
    has explicit error bounds and guarantees stability. However, it requires the
    solution of two Lyapunov equations to construct bases, which is a formidable
    task in large-scale problems. (ii) The moment-matching methods
    \cite{bai2002krylov, gugercin2008h_2} provide a computationally efficient
    framework using Krylov subspace techniques in an iterative fashion where
    only matrix-vector multiplications are required. The optimal $H_2$
    tangential interpolation for nonparametric systems \cite{gugercin2008h_2} is
    also available. The most crucial part of the moment-matching methods is
    location of samples where moments are matched. Also, it is not a data-driven
    approach, meaning that no data is used to construct ROM.  (iii) Proper
    Generalized Decomposition (PGD) \cite{ammar2006new} was first developed as a
    numerical method of solving boundary value problems, later extended to
    dynamical problems \cite{ammar2007new}. The main assumption of the method is
    a separated solution representation in space and time, which gives a way for
    an efficient solution procedure.  Therefore, it is considered as a model
    reduction technique.  However, PGD is not a data-driven approach. 

    Often there are many data available either from experiments or high-fidelity
    simulations. Those data contain valuable information about the system of
    interest.  Therefore, data-driven methods can maximize the usage of the
    existing data and enable the construction of an optimal ROM.  Dynamic Mode
    Decomposition (DMD) is a data-driven approach that generates reduced modes
    that embed an intrinsic temporal behavior. It was first developed by Peter
    Schmid in \cite{schmid2010dynamic} and populated by many other scientists.
    For more detailed information about DMD, we refer to this preprint
    \cite{tu2013dynamic}.
    As another data-driven approach, Proper Orthogonal decomposition (POD)
    \cite{berkooz1993proper} gives the data-driven optimal basis through the
    method of snapshots.  However, most of POD-based ROMs for linear dynamical
    systems apply spatial projection only.  Temporal complexity is still
    proportional to temporal discretization of its high-fidelity model. In order
    to achieve an optimal reduction, a space--time ROM needs to be built where
    both spatial and temporal projections are applied.  In literature, some
    space--time ROMs are available \cite{choi2019space, urban2014improved,
    yano2014space1, yano2014space2}, but they are only applied to small-scale
    problems. 


    Recently, several ROM techniques have been applied to various types of
    transport equations. Starting with Wols's work that used the POD to simulate
    the dynamics of an accelerator driven system (ADS) in 2010
    \cite{wols2010transient}, the publications on ROMs in transport problems
    have increased in number. For example, a POD-based ROM for eigenvalue
    problems to calculate dominant eigenvalues in reactor physics applications
    is developed by Buchan, et al. in \cite{buchan2013pod}.  The corresponding
    Boltzmann transport equation was re-casted into its diffusion form, where
    some of dimensions were eliminated.  Sartori, et al. in
    \cite{sartori2014comparison} also applied the POD-based ROM to the diffusion
    form and compared it with a modal method to show that the POD was superior
    to the modal method.  Reed and Robert in \cite{reed2015energy} also used POD
    to expand energy that replaced the traditional Discrete Legendre Polynomials
    (DLP) or modified DLP. They showed that a small number of POD energy modes
    could capture many-group fidelity.  Buchan, et al. in \cite{buchan2015pod}
    also developed a POD-based reduced order model to efficiently resolve the
    angular dimension of the steady-state, mono-energetic Boltzmann transport
    equation.  Behne, Ragusa, and Morel \cite{behne2019model} applied POD-based
    reduced order model to accelerate steady state $S_n$ radiation multi-group
    energy transport problem. A Petrov-Galerkin projection was used to close the
    reduced system.  Coale and Anistratov in \cite{coale2019areduced} replaces a
    high-order (HO) system with POD-based ROM in high-order low-order (HOLO)
    approach for Thermal Radiative Transfer (TRT) problems. 
   
    There have been some interesting DMD works for transport problems.
    McClarren and Haut in \cite{mcclarren2018acceleration} used DMD to estimate
    the slowly decaying modes from Richardson iteration and remove them from the
    solution.  Hardy, Morel, and Cory \cite{hardy2019dynamic} also explored DMD
    to accelerate the kinetics of subcritical metal systems without losing much
    accuracy. It was applied to a three-group diffusion model in a bare
    homogeneous fissioning sphere.  An interesting work by Star, et al. in
    \cite{star2019pod} exists, using the DMD method to identify non-intrusive
    POD-based ROM for the unsteady convection-diffusion scalar transport
    equation.  Their approach applied the DMD method to reduced coordinate
    systems.

    Several papers are found to use the PGD for transport problems. For example,
    Prince and Ragusa \cite{prince2019separated} applied the Proper Generalized
    Decomposition (PGD) to steady-state mono-energetic neutron transport
    equations where $S_n$ angular flux was sought as a finite sum of separable
    one-dimensional functions. However, the PGD was found to be ineffective for
    pure absorption problems because a large number of terms were required in
    the separated representation.  Prince and Ragusa \cite{prince2019parametric}
    also applied the PGD for uncertainty quantification process in the neutron
    diffusion--reaction problems.  Dominesey and Ji used the PGD method to
    separate space and angle in \cite{dominesey2019reduced} and to separate
    space and energy in \cite{dominesey2019areduced}.

    However, all these model order reduction techniques for transport equations
    apply only spatial projection, ignoring the potential reduction in temporal
    dimension.  In this paper, a Space--Time Reduced Order Model (ST-ROM) is
    developed and applied to large-scale linear dynamical problems.  Our ST-ROM
    achieves complexity reduction in both space and time dimension, which
    enables a great maximal speed-up and accuracy. It is amenable to any time
    integrators. It  follows the framework initially published in
    \cite{choi2019space}, but makes a new contribution by discovering a block
    structure in space--time basis that enables efficient implementation of the
    ST-ROM.  The block structure in space--time basis allows us not to build a
    space--time basis explicitly. It enables the construction of space--time
    reduced operators with small additional costs to the space-only reduced
    operators. In turn, this allows us to apply the space--time ROM to a
    large-scale linear dynamical problem. 

    The paper is organized in the following way: Section~\ref{sec:notations}
    introduces useful mathematical notations that will be used throughout the
    paper. Section~\ref{sec:lineardynamicalsystems} describes a parametric
    linear dynamical system and how to solve high-fidelity model in a classical
    time marching fashion. The full-order space--time formulation is also
    presented in Section~\ref{sec:lineardynamicalsystems} to be reduced to form
    our space--time ROM in Section~\ref{sec:spacetimeROM}. Section~\ref{sec:ROM}
    introduces both spatial and spatiotemporal ROMs. The basis generation is
    described in Section~\ref{sec:basis} where the traditional POD and
    incremental POD are explained in Sections~\ref{sec:POD} and
    \ref{sec:incremental}, respectively. Section~\ref{sec:blockstructure}
    reveals a block structure of the space--time reduced basis and derive each
    space--time reduced operators in terms of the blocks.  We apply our
    space--time ROM to a large-scale linear dynamical problem, i.e., a neutron
    transport simulation of solving BTE.  Section~\ref{sec:NeutronTransport}
    explains a discretization derivation of the Boltzmann transport equation,
    using multigroup energy discretization, associated Legendre polynomials for
    surface harmonic, simple corner balance discretization for space and
    direction, and the discrete ordinates method. Finally, we present our
    numerical results in Section~\ref{sec:numericalresults} and conclude the
    paper with summary and future works in Section~\ref{sec:conclusion}.

\subsection{Notations}\label{sec:notations}
  We review some of the notation used throughout the paper. 
  An $\ell_2$ norm is denoted as $\|\cdot\|$.
  For matrices
  $\boldsymbol{A} \in \RR{m \times n}$ and $\boldsymbol{B} \in \RR{k \times
  l}$, the {\em Kronecker} (or {\em tensor}) {\em product} of $\boldsymbol{A}$
  and $\boldsymbol{B}$ is the
  $mk \times nl$ matrix denoted by 
  \[
    \boldsymbol{A} \tensor \boldsymbol{B} \equiv 
  \left( \begin{array}{ccc}
    a_{11} \boldsymbol{B} & \cdots & a_{1n} \boldsymbol{B} \\
  \vdots & \ddots & \vdots \\
    a_{m1} \boldsymbol{B} & \cdots & a_{mn} \boldsymbol{B}
  \end{array} \right) ,
  \]
  where $\boldsymbol{A} = (a_{ij})$.  Kronecker products have many interesting
  properties.  We list here the ones relevant to our discussion:
  \begin{itemize}

    \item If $\boldsymbol{A}$ and $\boldsymbol{B}$ are nonsingular, then
      $\boldsymbol{A} \tensor \boldsymbol{B}$ is
        nonsingular with $(\boldsymbol{A} \tensor \boldsymbol{B})^{-1} =
        \boldsymbol{A}^{-1} \tensor \boldsymbol{B}^{-1}$,

      \item $(\boldsymbol{A} \tensor \boldsymbol{B})^T = \boldsymbol{A}^T
        \tensor \boldsymbol{B}^T$,

      \item Given matrices $\boldsymbol{A}$, $\boldsymbol{B}$, $\boldsymbol{C}$,
        and $\boldsymbol{D}$, $(\boldsymbol{A} \tensor \boldsymbol{B}) \cdot
        (\boldsymbol{C}
        \tensor \boldsymbol{D}) = \boldsymbol{AC} \tensor \boldsymbol{BD}$, as long as both sides of the
        equation make sense,

      \item $(\boldsymbol{A} + \boldsymbol{B}) \tensor \boldsymbol{C} =
        \boldsymbol{A} \tensor \boldsymbol{C} + \boldsymbol{B} \tensor
        \boldsymbol{C}$, and

      \item $\boldsymbol{A} \tensor (\boldsymbol{B}+\boldsymbol{C}) =
        \boldsymbol{A} \tensor \boldsymbol{B} + \boldsymbol{A} \tensor
        \boldsymbol{C}$.
  \end{itemize}

\section{Linear dynamical systems}\label{sec:lineardynamicalsystems}
  Parametric continuous dynamical systems that are linear in state 
  are considered:
  \begin{align}
    \dot{\state}(\timevar;\param) &= \statemat(\param) \state(\timevar;\param) +
    \inputmat(\param)\inputvec(\timevar;\param), \quad 
    \state(0;\param) = \initialstate(\param),\label{eq:Dynamicalstate} \\  
    \outputvec(\timevar;\param) &= \outputmat(\param)^T \state(\timevar;\param), 
    \label{eq:Dynamicaloutput}
  \end{align}
  where $\param\in\paramspace\subset\RR{\nparam}$ denotes a parameter vector,
  $\state: [0,\finaltime]\times\RR{\nparam} \rightarrow \RR{\nspace}$ denotes a
  time dependent state variable function,
  $\initialstate:\RR{\nparam}\rightarrow\RR{\nspace}$ denotes an initial state,
  $\inputvec: [\initialtime,\finaltime]\times\RR{\nparam} \rightarrow
  \RR{\ninput}$ denotes a time dependent input variable function, and
  $\outputvec: [\initialtime, \finaltime]\times\RR{\nparam} \rightarrow
  \RR{\noutput}$ denotes a time dependent output variable function. The system
  operations, i.e.,
  $\statemat:\RR{\nparam}\rightarrow\RR{\nspace\times\nspace}$,
  $\inputmat:\RR{\nparam}\rightarrow\RR{\nspace\times\ninput}$, and
  $\outputmat:\RR{\nparam}\rightarrow\RR{\nspace\times\noutput}$, are real
  valued matrices, independent of state variables.  We assume that the dynamical
  system above is stable, i.e., the eigenvalues of $\statemat$ have strictly
  negative real parts.

  Our methodology works for any time integrators, but for the illustration
  purpose, we apply a backward Euler time integrator to Eq.~\eqref{eq:Dynamicalstate}.
  At $k$th time step, the following system of equations is solved:
  \begin{equation}\label{eq:fom}
    \left (\identityk{\nspace} - \timestepk{k}\statemat(\param) \right )\statek{k} = \statek{k-1} +
    \timestepk{k}\inputmat(\param)\inputveck{k}(\param),
  \end{equation}
  where $\identityk{\nspace}\in\RR{\nspace\times\nspace}$ denotes an identity
  matrix, $\timestepk{k}$ denotes $k$th time step size with $\finaltime =
  \sum_{k=1}^{\ntime} \timestepk{k}$, and $\statek{k}(\param)$ and
  $\inputveck{k}(\param)$
  denote state and input vectors at $k$th time step,
  $\timevark{k} = \sum_{j=1}^k \timestepk{j}$, respectively.   
  A Full Order Model (FOM) solves Eq.~\eqref{eq:fom} every time step. 
  The spatial dimension, $\nspace$, and the temporal dimension, $\ntime$ can be
  very large, which leads to a large-scale problem. We introduce how to reduce
  the high dimensionality in Section~\ref{sec:ROM}.

  The single time step formulation in Eq.~\eqref{eq:fom} can be equivalently
  re-written in the following discretized space-time formulation:
  \begin{equation}\label{eq:st-fom}
    \ststatemat(\param)\ststate(\param) = \stinputvec(\param) +
    \stinitial(\param),
  \end{equation}
  where the space–-time system matrix, $\ststatemat:\RR{\nparam} \rightarrow
  \RR{\nspace\ntime \times \nspace\ntime}$, the space--time state vector,
  $\ststate:\RR{\nparam}\rightarrow\RR{\nspace\ntime}$, the space--time input
  vector, $\stinputvec:\RR{\nparam}\rightarrow\RR{\nspace\ntime}$, and the
  space--time initial vector,
  $\stinitial:\RR{\nparam}\rightarrow\RR{\nspace\ntime}$, are defined
  respectively as 
  \begin{equation}\label{eq:stsystemmat}
    \ststatemat(\param) = \bmat{
      \identityk{\nspace} - \timestepk{1} \statemat(\param) &  &    
    \\ -\identityk{\nspace} & \identityk{\nspace} - \timestepk{2}
    \statemat(\param) & & 
\\    & \ddots & \ddots & 
    \\   &   & -\identityk{\nspace} & \identityk{\nspace} -
    \timestepk{\ntime} \statemat(\param)
    },
  \end{equation}
  \begin{equation}\label{eq:ststate}
    \ststate(\param)  = \bmat{ \statek{1}(\param) \\ \statek{2}(\param) \\
    \vdots \\ \statek{\ntime}(\param) 
    },\quad
    \stinputvec(\param) = \bmat{
      \timestepk{1} \inputmat(\param)\inputveck{1}(\param)  \\ \timestepk{2}
      \inputmat(\param)\inputveck{2}(\param) \\ \vdots \\ \timestepk{\ntime}
      \inputmat(\param)\inputvec{\ntime}(\param) 
    },\quad
    \stinitial(\param) = \bmat{
      \initialstate(\param) \\ \zero \\ \vdots \\ \zero
    }.
  \end{equation}
  A lower block-triangular matrix structure of $\ststatemat$ comes from the
  backward Euler time integration scheme. Other time integrators will give other
  sparse block structures. No one will solve this space--time system directly
  because the specific block structure of $\ststatemat$ lets one to solve the
  system in time-marching fashion. However, if the space--time formulation in
  Eq.~\eqref{eq:st-fom} can be reduced and solved efficiently, then one might be
  interested in solving the reduced space--time system in its whole.
  Section~\ref{sec:spacetimeROM} shows such a reduction is possible. 

\section{Reduced order models}\label{sec:ROM}
  We consider a projection-based reduced order model for linear dynamical
  systems. Section~\ref{sec:spaceROM} shows a typical spatial reduced order
  model.  A space--time reduced order model is described in
  Section~\ref{sec:spacetimeROM}.

\subsection{Spatial reduced order models}\label{sec:spaceROM}
  A projection-based spatial reduced order model approximates the state
  variables as a linear combination of a small number of spatial basis vectors,
  $\{ \spatialbasisveck{1}, \ldots, \spatialbasisveck{\nreducedspace} \}$, where
  $\spatialbasisveck{k}\in\RR{\nspace}$, $k\in\nat{\nreducedspace}$
  with $\nat{\nreducedspace} \equiv \{1,\ldots,\nreducedspace\}$,
  $\nreducedspace \ll \nspace$, i.e.,
  \begin{equation}\label{eq:stateapproximation}
    \state(\timevar;\param) \approx \approxstate(\timevar) \equiv \refstate(\param) +
    \spatialbasis\reducedstate(\timevar;\param),
  \end{equation}
  where the spatial basis, $\spatialbasis\in\RR{\nspace*\nreducedspace}$ with
  $\spatialbasis^T\spatialbasis = \identityk{\nreducedspace}$, is
  defined as 
  \begin{equation}\label{eq:sptialbasis}
    \spatialbasis \equiv \bmat{\spatialbasisveck{1} & \cdots &
    \spatialbasisveck{\nreducedspace}},
  \end{equation}
  a reference state is denoted as $\refstate(\param)\in\RR{\nspace}$, and a
  time-dependent reduced coordinate vector function is defined as
  $\reducedstate: \RR{\nparam} \rightarrow \RR{\nreducedspace}$.
  Substituting \eqref{eq:stateapproximation} to \eqref{eq:Dynamicalstate}, gives an
  over-determined system of equations:
  \begin{equation}\label{eq:overdeterminedODE}
    \spatialbasis\dot{\reducedstate}(\timevar;\param) =
    \statemat(\param)\refstate(\param) +
    \statemat(\param)
    \spatialbasis\reducedstate(\timevar;\param) +
    \inputmat(\param)\inputvec(\timevar,\param), 
  \end{equation}
  which can be closed, for example, by Galerkin projection, i.e.,
  left-multiplying both sides of \eqref{eq:overdeterminedODE} and initial
  condition in \eqref{eq:Dynamicalstate} by $\spatialbasis^T$, giving the reduced
  system of equations and initial conditions:
  \begin{equation}\label{eq:galerkin}
    \dot{\reducedstate}(\timevar;\param) = \reducedstatemat(\param)
    \reducedstate(\timevar;\param) +
    \reducedinputmat(\param)\inputvec(\timevar;\param) +
    \reducedrefstate(\param), \quad \reducedstate(0;\param) =
    \reducedinitialstate(\param)
  \end{equation}
  where $\reducedstatemat(\param) \equiv \spatialbasis^T\statemat(\param)
  \spatialbasis \in \RR{\nreducedspace\times\nreducedspace}$ denotes a reduced
  system matrix, $\reducedinputmat(\param) \equiv \spatialbasis^T\inputmat(\param) \in
  \RR{\nreducedspace\times\ninput}$ denotes a reduced input matrix,
  $\reducedrefstate(\param) \equiv \spatialbasis^T\statemat(\param)\refstate(\param) \in
  \RR{\nreducedspace}$ denotes a reduced reference state vector, and
  $\reducedinitialstate(\param) \equiv
  \spatialbasis^T(\initialstate(\param)-\refstate(\param))\in\RR{\nreducedspace}$
  denotes a reduced initial condition. Once $\spatialbasis$ and $\refstate$ are
  known, the reduced operators, $\reducedstatemat$, $\reducedinputmat$,
  $\reducedrefstate$, and $\reducedinitialstate$ can be pre-computed. With the
  pre-computed operators, the system~\eqref{eq:galerkin} can be solved fairly
  quickly, for example, by applying the backward Euler time integrator:
  \begin{equation}\label{eq:space_reduced_discretized_ODE}
    \left (\identityk{\nreducedspace}
    -\timestepk{k} \reducedstatemat(\param) \right ) \reducedstatek{k}(\param) 
    = \reducedstatek{k-1}(\param) + \timestepk{k}\reducedinputmat(\param)\inputveck{k}(\param) 
    + \timestepk{k}\reducedrefstate(\param), 
  \end{equation}
  where $\reducedstatek{k}(\param) \equiv
  \reducedstate(\timevark{k};\param)\in\RR{\nreducedspace}$.  Then the output
  vector, $\outputveck{k}(\param) \equiv \outputvec(\timevark{k};\param)$, can be
  computed as 
  \begin{equation}\label{eq:outputFromReducedsol}
    \outputveck{k}(\param) \equiv
    \reducedoutputmat(\param)^T\reducedstatek{k}(\param) +
    \outputmat(\param)^T\refstate(\param),
  \end{equation}
  where  $\reducedoutputmat(\param) \equiv \spatialbasis^T\outputmat(\param) \in
  \RR{\nreducedspace\times\noutput}$ denotes a reduced output matrix. 
  The usual choices for $\refstate$ include $\zerobold$, $\initialstate$, and
  some kind of average quantities. Note that if $\initialstate$ is used as
  $\refstate$, then $\reducedinitialstate = \zerobold$, independent of $\param$,
  which is convenient.  The POD for generating the spatial basis is described in
  Section~\ref{sec:basis}.

\subsection{Space--time reduced order models}\label{sec:spacetimeROM}
  The space--time formulation, \eqref{eq:st-fom}, can be reduced by
  approximating the space--time state variables as a linear combination of a
  small number of space--time basis vectors,  
  $\{ \spacetimebasisveck{1}, \ldots,
  \spacetimebasisveck{\nreducedspace\nreducedtime} \}$, where
  $\spacetimebasisveck{k}\in\RR{\nspace\ntime}$,
  $k\in\nat{\nreducedspace\nreducedtime}$
  with $\nreducedspace\nreducedtime \ll \nspace\ntime$, i.e.,
  \begin{equation}\label{eq:ststateapproximation}
    \ststate(\param) \approx \approxststate(\param) \equiv 
    \spacetimebasis\reducedststate(\param),
  \end{equation}
  where the space--time basis,
  $\spacetimebasis\in\RR{\nspace\ntime*\nreducedspace\nreducedtime}$ is
  defined as 
  \begin{equation}\label{eq:spacetimebasis}
    \spacetimebasis \equiv \bmat{\spacetimebasisveck{1} & \cdots &
    \spacetimebasisveck{\reducedspaceindex+\nreducedspace(\reducedtimeindex-1)}  & \cdots
    \spacetimebasisveck{\nreducedspace\nreducedtime}},
  \end{equation}
  where $\reducedspaceindex\in\nat{\nreducedspace}$,
  $\reducedtimeindex\in\nat{\nreducedtime}$.  The space--time reduced coordinate
  vector function is denoted as $\reducedststate:\RR{\nparam} \rightarrow
  \RR{\nreducedspace\nreducedtime}$.  Substituting
  \eqref{eq:ststateapproximation} to \eqref{eq:st-fom}, gives an over-determined
  system of equations:
  \begin{equation}\label{eq:overdetermined-spacetime}
    \ststatemat(\param)\spacetimebasis\reducedststate(\param) = \stinputvec(\param) +
    \stinitial(\param) 
  \end{equation}
  which can be closed, for example, by Galerkin projection, i.e.,
  left-multiplying both sides of \eqref{eq:overdetermined-spacetime} 
  by $\spacetimebasis^T$, giving the reduced
  system of equations:
  \begin{equation}\label{eq:st-galerkin}
    \reducedststatemat(\param)\reducedststate(\param) =
    \reducedstinputvec(\param) + \reducedstinitial(\param) 
  \end{equation}
  where $\reducedststatemat(\param) \equiv \spacetimebasis^T\ststatemat(\param)
  \spacetimebasis \in
  \RR{\nreducedspace\nreducedtime\times\nreducedspace\nreducedtime}$ denotes a
  reduced space--time system matrix, $\reducedstinputvec(\param) \equiv
  \spacetimebasis^T\stinputvec(\param) \in \RR{\nreducedspace\nreducedtime}$ denotes a
  reduced space--time input vector, and $\reducedstinitial(\param) \equiv
  \spacetimebasis^T\stinitial(\param) \in \RR{\nreducedspace\nreducedtime}$ denotes a
  reduced space--time initial state vector.  Once $\spacetimebasis$ is known,
  the space--time reduced operators, $\reducedststatemat$,
  $\reducedstinputvec$, and $\reducedstinitial$, can be pre-computed, but its
  computation involves the
  number of operations in $O(\nspace\ntime)$, which can be large.
  Sec.~\ref{sec:blockstructure} explores a block structure of $\spacetimebasis$
  that shows an efficient way of constructing reduced space--time operators
  without explicitly forming the full-size space--time operators, such as
  $\spacetimebasis$, $\ststatemat$, $\stinputvec$, and $\stinitial$.

\section{Basis generation}\label{sec:basis}
\subsection{Proper orthogonal decomposition}\label{sec:POD}

 \begin{figure}[t!]
    \centering
    \includegraphics[width=0.8\textwidth]{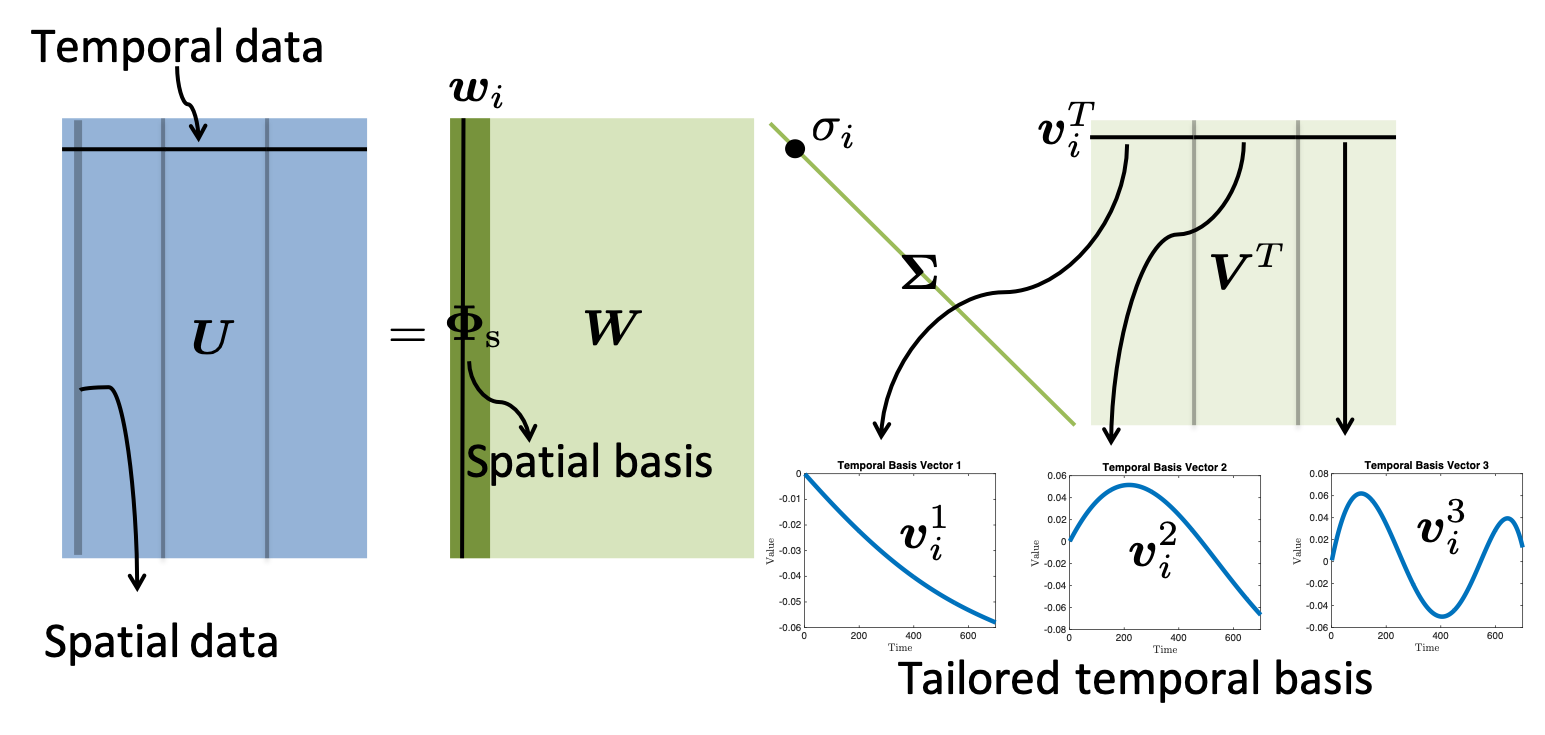}
    \caption{Illustration of spatial and temporal bases construction, using SVD
    with $\nparam=3$. The right singular vector,
    $\rightsingularveck{\reducedspaceindex}$,
    describes three different temporal behaviors of a left singular basis vector
    $\leftsingularveck{\reducedspaceindex}$, i.e., three different temporal behaviors of a
    spatial mode. Each temporal behavior is denoted as
    $\rightsingularveck{\reducedspaceindex}^1$,
    $\rightsingularveck{\reducedspaceindex}^2$, and
    $\rightsingularveck{\reducedspaceindex}^3$.  }
    \label{fig:SVDexplanation}
 \end{figure}

  We follow the method of snapshots first introduced by Sirovich
  \cite{sirovich1987turbulence}.  Let $\paramsample =
  \{\paramk{1},\ldots,\paramk{\nparam}\}$ be a set of parameter samples where we
  run full order model simulations. Let $\stateSnapshotk{p} \equiv
  \bmat{\statek{1}(\paramk{p}) & \cdots &
  \statek{\ntime}(\paramk{p})}\in\RR{\nspace\times\ntime}$,
  $p\in\nat{\nparam}$, be a full order model state solution matrix for a sample
  parameter value, $\paramk{p}\in\paramspace$. Then a snapshot matrix,
  $\stateSnapshot\in\RR{\nspace\times\nparam\ntime}$, is defined by
  concatenating all the state solution matrices, i.e., 
  \begin{equation}\label{eq:snapshotmat}
    \stateSnapshot \equiv \bmat{\stateSnapshotk{1} & \cdots &
    \stateSnapshotk{\nparam}}.
  \end{equation}
  The spatial basis from POD is an optimally compressed
  representation of $\range{\stateSnapshot}$ in a sense that it minimizes the
  difference between the original snapshot matrix and the projected one onto the
  subspace spanned by the basis, $\spatialbasis$:
  \begin{equation}\label{eq:POD}
    \begin{aligned}
      & \underset{\spatialbasis\in\RR{\nspace\times\nreducedspace},
      \spatialbasis^T\spatialbasis
      =\identityk{\nreducedspace}}{\text{minimize}} & & \left \|\stateSnapshot -
      \spatialbasis\spatialbasis^T{\stateSnapshot} \right \|_F^2, 
    \end{aligned}
  \end{equation}
  where $\|\cdot\|_F$ denotes the Frobenius norm.  The solution of POD can be
  obtained by setting $\spatialbasis = \leftsingularmat(:,1:\nreducedspace)$,
  $\nreducedspace < \nparam\ntime$, in
  MATLAB notation, where $\leftsingularmat$ is the left singular matrix of the
  following thin Singular Value Decomposition (SVD) with $\ell \equiv
  \min(\nspace,\nparam\ntime)$:
 \begin{align}
   \stateSnapshot &= \leftsingularmat\singularvalmat\rightsingularmat^T
   \label{eq:SVD}\\
                  &= \sum_{\reducedspaceindex=1}^{\ell}\singularvaluesArg{i}
                  \leftsingularveck{\reducedspaceindex}
                  \rightsingularveck{\reducedspaceindex}^T\label{eq:SVDsummation}
 \end{align} 
 where $\leftsingularmat\in\RR{\nspace\times\ell}$ and
 $\rightsingularmat\in\RR{\nparam\ntime\times\ell}$ are orthogonal matrices and
 $\singularvalmat\in\RR{\ell\times\ell}$ is a diagonal matrix with singular
 values on its diagonal. The equivalent summation form is written in
 \eqref{eq:SVDsummation}, where $\singularvaluesArg{\reducedspaceindex}\in\RR{}$
 is $\reducedspaceindex$th singular value,
 $\leftsingularveck{\reducedspaceindex}$ and $\rightsingularveck{i}$ are $i$th
 left and right singular vectors,
 respectively. Note that $\rightsingularveck{\reducedspaceindex}$ describes $\nparam$ different
 temporal behavior of $\leftsingularveck{\reducedspaceindex}$. For example,
 Figure~\ref{fig:SVDexplanation} illustrates the case of $\nparam=3$, where
 $\rightsingularveck{\reducedspaceindex}^1$,
 $\rightsingularveck{\reducedspaceindex}^2$, and
 $\rightsingularveck{\reducedspaceindex}^3$ describe three different temporal
 behavior of a specific spatial basis vector, i.e.,
 $\leftsingularveck{\reducedspaceindex}$. For general $\nparam$, we note that
 $\rightsingularveck{\reducedspaceindex}$ describes $\nparam$ different temporal
 behavior of $\reducedspaceindex$th spatial basis vector, i.e.,
 $\spatialbasisveck{\reducedspaceindex} =
 \leftsingularveck{\reducedspaceindex}$. We set
 $\temporalSnapshotk{\reducedspaceindex} =
 \bmat{\rightsingularveck{\reducedspaceindex}^1 & \cdots &
 \rightsingularveck{\reducedspaceindex}^{\nparam}   }$
 to be $\reducedspaceindex$th temporal snapshot matrix, where
 $\rightsingularveck{\reducedspaceindex}^k \equiv
 \rightsingularveck{\reducedspaceindex}(1+(k-1)\ntime:k\ntime)$ for $k\in\nat{\nparam}$.  We
 apply SVD on $\temporalSnapshotk{\reducedspaceindex}$:
 \begin{equation}\label{eq:temporalSVD}
   \temporalSnapshotk{\reducedspaceindex} =
   \leftsingularmattemporalk{\reducedspaceindex}\singularvalmatk{\reducedspaceindex}\rightsingularmattemporalk{\reducedspaceindex}^T.
 \end{equation}
 Then, the temporal basis for $\reducedspaceindex$th spatial basis vector can be
 set
 $\temporalbasisi{\reducedspaceindex} =
 \leftsingularmattemporalk{\reducedspaceindex}(:,1:\nreducedtime)$ in
 MATLAB notation. Finally, a space--time basis vector,
 $\spacetimebasisveck{\reducedspaceindex+\nreducedspace(\reducedtimeindex-1)}\in\RR{\nspace\ntime}$, in
 \eqref{eq:spacetimebasis} can be constructed as
  \begin{equation}\label{eq:spacetimebasisvec}
    \spacetimebasisveck{\reducedspaceindex+\nreducedspace(\reducedtimeindex-1)} =
    \temporalbasisveck{\reducedspaceindex}{\reducedtimeindex}\otimes\spatialbasisveck{\reducedspaceindex},
  \end{equation}
  where $\spatialbasisveck{\reducedspaceindex} \equiv
  \spatialbasis(:,\reducedspaceindex) \in \RR{\nspace}$ denotes
  $\reducedspaceindex$th spatial basis vector and
  $\temporalbasisveck{\reducedspaceindex}{\reducedtimeindex} \equiv
  \temporalbasisi{\reducedspaceindex}(:,\reducedtimeindex)\in\RR{\ntime}$
  denotes $\reducedtimeindex$th temporal basis vector that describes a temporal
  behavior of $\spatialbasisveck{\reducedspaceindex}$. The computational cost of
  SVD for the snapshot matrix,
  $\stateSnapshot\in\RR{\nspace\times\nparam\ntime}$, assuming $\nspace \gg
  \nparam\ntime$, is $O(\nspace^2\nparam\ntime)$ and the computational cost of
  SVD for $\nreducedspace$ temporal snapshot matrices,
  $\temporalSnapshotk{i}\in\RR{\ntime\times\nparam}$, $i\in\nat\nreducedspace$,
  $\ntime \gg \nparam$ is $O(\ntime^2\nreducedspace\nparam)$. For a large-scale
  problem, this may be a formidable task.  Thus, we use an incremental SVD where
  a rank one update of existing SVD is achieved with much more memory-efficient
  way than the thin SVD in Eq.~\eqref{eq:SVD}.  The incremental SVD procedure is
  explained in Section~\ref{sec:incremental}.

 POD is related to the principal component analysis in
 statistics \cite{hotelling1933analysis} and Karhunen--Lo\`{e}ve expansion
 \cite{loeve1955} in stochastic analysis.  Since the objective function in
 \eqref{eq:POD} does not change even though $\spatialbasis$ is post-multiplied
 by an arbitrary $\nreducedspace\times\nreducedspace$ orthogonal matrix, the POD
 procedure seeks the optimal $\nreducedspace$-dimensional subspace that captures
 the snapshots in the least-squares sense.  For more details on POD, we refer to
 \cite{berkooz1993proper,hinze2005proper,kunisch2002galerkin}.

\subsection{Incremental space--time reduced basis}\label{sec:incremental}
 An incremental SVD is an efficient way of updating the existing singular value
 decomposition when a new snapshot vector, i.e., a column vector, is added. For
 a time dependent problem, we start with a first time step solution with a first
 parameter vector, i.e., $\statek{1}(\paramk{1})$.  If its norm is big enough
 (i.e., $\|\statek{1}(\paramk{1})\| > \svdtol$), then we set the first singular
 value $\singularvaluesArg{1} = \|\statek{1}(\paramk{1})\|$, the first
 left singular vector be the normalized first snapshot vector, i.e.,
 $\leftsingularveck{1} = \statek{1}(\paramk{1})/\singularvaluesArg{1}$, and the
 right singular vector be $\rightsingularveck{1} = 1$. Otherwise, we set them
 empty, i.e., $\singularvaluesArg{1} = []$, $\leftsingularveck{1} = []$, and
 $\rightsingularveck{1} = []$.  This initializing process is described in
 Algorithm~\ref{al:initializingIncrementalSVD}. We pass $k$ to
 \textbf{initializingIncrementalSVD} function as an input argument to indicate
 $k$th snapshot vector is being handled. Also, the rank of
 $\leftsingularmatk{k}$ is denoted as $\basisrankk{k}$. In general,
 $\basisrankk{k} \neq k$ because a snapshot vector will not be included if
 it is too small (i.e., Line 1 in Algorithm~\ref{al:initializingIncrementalSVD})
 or it is linearly dependent on the existing basis (i.e., Line 9 and 13 in
 Algorithm~\ref{al:incrementalSVD}) or it generates a small eigenvalue (i.e.,
 Line 18 in Algorithm~\ref{al:incrementalSVD}).
 
 Let's assume that we have $(k-1)$th SVD from previous $k-1$ snapshot vectors,
 i.e., $\leftsingularmatk{k-1} \singularvalmat_{k-1} \rightsingularmatk{k-1}^T$,
 whose rank is $\basisrankk{k-1}$.
 If a new snapshot vector, $\snapshotvec$ (e.g., $k$th time step solution with the first sample
 parameter value, $\statek{k}(\paramk{1})$) needs to be added to the existing
 SVD, the following factorization can be used \cite{brand2002incremental}:
  \begin{align}
    \bmat{ \leftsingularmatk{k-1} \singularvalmat_{k-1}
    \rightsingularmatk{k-1}^T & \snapshotvec}
    &= \bmat{ \leftsingularmatk{k-1} & \left
    (\identity-\leftsingularmatk{k-1}\leftsingularmatk{k-1}^T\right)
    \snapshotvec/\projerr} 
    \bmat{\singularvalmat_{k-1} &
    \leftsingularmatk{k-1}^T\snapshotvec \\ \zerobold & \projerr}
    \bmat{ \rightsingularmatk{k-1} & \zerobold \\ \zerobold & 1}^T
    \label{eq:factorization1}\\
    &= \bmat{ \leftsingularmatk{k-1} & \updatedBasisVector }
    \bmat{\singularvalmat_{k-1} & \reducedsnapshotvec \\ \zerobold
    & \projerr}
    \bmat{ \rightsingularmatk{k-1} & \zerobold \\ \zerobold & 1}^T,
    \label{eq:factorization2}
  \end{align}
  where $\reducedsnapshotvec = \leftsingularmatk{k-1}^T\snapshotvec \in
  \RR{\basisrankk{k-1}}$ denotes a
  reduced coordinate of $\snapshotvec$ that is projected onto the subspace
  spanned by $\leftsingularmatk{k-1}$, $\projerr = \| \snapshotvec -
  \leftsingularmatk{k-1}\reducedsnapshotvec \|$ denotes the norm of the
  difference between $\snapshotvec$ and the projected one, and
  $\updatedBasisVector = \left
  (\snapshotvec-\leftsingularmatk{k-1}\reducedsnapshotvec\right)/\projerr \in
  \RR{\nspace}$
  denotes a new orthogonal vector due to the incoming vector, $\snapshotvec$.
  Note that the left and right matrices of the factorization, i.e., $\bmat{
    \leftsingularmatk{k-1} & \updatedBasisVector } \in
  \RR{\nspace\times(\basisrankk{k-1}+1)}$ and  $\bmat{
    \rightsingularmatk{k-1} & \zerobold \\ \zerobold & 1} \in
  \RR{k\times(\basisrankk{k-1}+1)}$ are orthogonal
  matrices.  Let
  $\middlemat\in{\RR{(\basisrank_{k-1}+1)\times(\basisrank_{k-1}+1)}}$ denote
  the middle matrix of the factorization, i.e., 
  \begin{equation}\label{eq:middleMat}
    \middlemat = \bmat{\singularvalmat_{k-1} & \reducedsnapshotvec \\ \zerobold &
    \projerr}.
  \end{equation}
  The matrix, $\middlemat$, is almost diagonal except for $\reducedsnapshotvec$
  in the upper right block, i.e., one column bordered diagonal. Its size is not
  in $O(\nspace)$. Thus, the SVD of $\middlemat$ is computationally cheap, i.e.,
  $O((\basisrank_{k-1}+1)^3)$. Let the SVD of $\middlemat$ be
  \begin{equation}\label{eq:svdMiddleMat}
    \middlemat = \middleleftsingularmat\middlesingularvalmat
    \middlerightsingularmat, 
  \end{equation}
  where $\middleleftsingularmat \in \RR{(\basisrank_{k-1}+1)\times(\basisrank_{k-1}+1)}$
  denotes the left singular matrix, $\middlesingularvalmat \in
  \RR{(\basisrank_{k-1}+1)\times(\basisrank_{k-1}+1)}$ denotes the singular
  value matrix, and $\middlerightsingularmat \in
  \RR{(\basisrank_{k-1}+1)\times(\basisrank_{k-1}+1)}$ denotes the right
  singular matrix of $\middlemat$. Replacing $\middlemat$ in
  Eq.~\eqref{eq:factorization2} with \eqref{eq:svdMiddleMat} gives
  \begin{align}
    \bmat{ \leftsingularmatk{k-1} \singularvalmat_{k-1}
    \rightsingularmatk{k-1}^T & \snapshotvec} &= 
    \bmat{ \leftsingularmatk{k-1} & \updatedBasisVector }
    \middleleftsingularmat\middlesingularvalmat
    \middlerightsingularmat
    \bmat{ \rightsingularmatk{k-1} & \zerobold \\ \zerobold & 1}^T 
    \label{eq:rankOneUpdateSVD1}\\
    &= \leftsingularmatk{k} \singularvalmat_{k}
    \rightsingularmatk{k}^T \label{eq:rankOneUpdateSVD2},
  \end{align}
  where $\leftsingularmatk{k} = \bmat{ \leftsingularmatk{k-1} &
  \updatedBasisVector } \middleleftsingularmat \in
  \RR{\nspace\times(\basisrank_{k-1}+1)}$ denotes the updated left singular matrix,
  $\singularvalmat_{k} = \middlesingularvalmat \in
  \RR{(\basisrank_{k-1}+1)\times(\basisrank_{k-1}+1)}$ denotes the updated singular value
  matrix, and $\rightsingularmatk{k} = \bmat{ \rightsingularmatk{k-1} &
  \zerobold \\ \zerobold & 1}\middlerightsingularmat \in
  \RR{k\times(\basisrank_{k-1}+1)}$ denotes the updated right singular
  matrix. This updating algorithm is described in
  Algorithm~\ref{al:incrementalSVD}.

  Algorithm~\ref{al:incrementalSVD} also checks if $\snapshotvec$ is linearly
  dependent on the current left singular vectors numerically. If $\projerr <
  \svdtol$, then we consider that it is linearly dependent. Thus, we set
  $\projerr=0$ in $\middlemat$, i.e., Line 10 of
  Algorithm~\ref{al:incrementalSVD}.  Then we only update the first
  $\basisrank_{k-1}$ components of the singular matrices in Line 14 of
  Algorithm~\ref{al:incrementalSVD}. Otherwise, we follow the update form in
  Eq.~\eqref{eq:rankOneUpdateSVD2} as in Line 16 in
  Algorithm~\ref{al:incrementalSVD}.  
  
  Line 18-20 in Algorithm~\ref{al:incrementalSVD} checks if the updated singular
  value has a small value. If it does, we neglect that particular singular value
  and corresponding component in left and right singular matrices. It is because
  a small singular value causes a large error in left and right singular
  matrices \cite{fareed2018error}.  
  
  Although the orthogonality of the updated left singular matrix,
  $\leftsingularmatk{k}$, must be guaranteed in infinite precision by the
  product of two orthogonal matrices in Line 14 or 16 of
  Algorithm~\ref{al:incrementalSVD}, it is not guaranteed in finite precision.
  Thus, we heuristically check the orthogonality in Lines 21-24 of
  Algorithm~\ref{al:incrementalSVD} by checking the inner product of the first
  and last columns of $\basis_k$. If the orthogonality is not shown, then we
  orthogonalize them by the QR factorization. Here $\tol$ denotes unit roundoff
  (e.g., \textsf{eps} in MATLAB).
  
  The spatial basis can be set after $\nparam\ntime$ incremental steps:
  \begin{equation}\label{eq:spatialincrementalbasis}
    \spatialbasis = \leftsingularmatk{\nparam\ntime}(:,1:\nreducedspace).
  \end{equation}
  If all the time step solutions are taken incrementally and sequentially from
  $\nparam$ different high-fidelity time dependent simulations, then the right
  singular matrix, $\rightsingularmatk{\ntime\nparam} \in
  \RR{\ntime\nparam\times\basisrankk{\ntime\nparam}}$, holds $\nparam$ different
  temporal behavior for each spatial basis vector. For example,
  $\rightsingularveck{\reducedspaceindex}$ describes $\nparam$ different
  temporal behavior of $\leftsingularveck{\reducedspaceindex}$.  As in
  Section~\ref{sec:POD}, $\reducedspaceindex$th temporal snapshot matrix
  $\temporalSnapshotk{\reducedspaceindex}\in\RR{\ntime\times\nparam}$ can be
  defined as
  \begin{equation}\label{eq:incrementaltemporalsnapshotmat}
    \temporalSnapshotk{\reducedspaceindex} \equiv \bmat{
      \rightsingularveck{\reducedspaceindex}^1 & \cdots &
         \rightsingularveck{\reducedspaceindex}^{\nparam} },
  \end{equation}
  where $\rightsingularveck{\reducedspaceindex}^k \equiv
  \rightsingularmatk{\ntime\nparam}(1+(k-1)\ntime:k\ntime,\reducedspaceindex)$ for
  $k\in\nat{\nparam}$. If we take the SVD of
  $\temporalSnapshotk{\reducedspaceindex} =
  \leftsingularmattemporalk{\reducedspaceindex}
  \singularvalmatk{\reducedspaceindex}
  \rightsingularmattemporalk{\reducedspaceindex}^T$,
  then  
  the temporal basis for $\reducedspaceindex$th spatial basis vector can be set 
  \begin{equation}\label{eq:temporalincrementalbasis}
    \temporalbasisi{\reducedspaceindex} =
    \leftsingularmattemporalk{\reducedspaceindex}(:,1:\nreducedtime).
  \end{equation}

  \begin{algorithm}[t!]
    \caption{Initializing incremental SVD}\label{al:initializingIncrementalSVD}
    [$\leftsingularmatk{k}$, $\singularvaluesk{k}$, $\rightsingularmatk{k}$] =
    \textbf{initializingIncrementalSVD}($\snapshotvec$, $\svdtol$, $k$) \\
    \textbf{Input:} $\snapshotvec$, $\svdtol$, $k$\\
    \textbf{Output:} $\leftsingularmatk{k}$, $\singularvaluesk{k}$,
    $\rightsingularmatk{k}$
      \begin{algorithmic}[1]
        \IF{$\|\snapshotvec\| > \svdtol$}
          \STATE $\singularvaluesk{k} \assign  \bmat{\|\snapshotvec\|}$,
          $\leftsingularmatk{k} \assign \snapshotvec/ \singularvaluesArg{1}$, and
          $\rightsingularmatk{k} \assign \bmat{1}$
        \ELSE
          \STATE $\singularvaluesk{k} \assign []$,
          $\leftsingularmatk{k} \assign []$, and 
          $\rightsingularmatk{k} \assign []$
        \ENDIF
   \end{algorithmic}
  \end{algorithm}

  \begin{algorithm}[t!]
    \caption{Incremental SVD, $\leftsingularmatk{-1} = [ ]$}\label{al:incrementalSVD}
    [$\leftsingularmatk{k}$, $\singularvaluesk{k}$, $\rightsingularmatk{k}$] =
    \textbf{incrementalSVD}($\snapshotvec$, $\svdtol$, $\svtol$, $\leftsingularmatk{k-1}$,
    $\singularvaluesk{k}$, $\rightsingularmatk{k-1}$, $k$)\\
    \textbf{Input:} $\snapshotvec$, $\svdtol$, $\svtol$, $\leftsingularmatk{k-1}$,
    $\singularvaluesk{k-1}$, $\rightsingularmatk{k-1}$, $k$\\
    \textbf{Output:} $\leftsingularmatk{k}$, $\singularvaluesk{k}$,
    $\rightsingularmatk{k}$
      \begin{algorithmic}[1]
        \IF{$\basisrank_{k-1}=0$ or $\basisrank_{k-1}=\maxrank$}
          \STATE [$\leftsingularmatk{k}$, $\singularvaluesk{k}$,
          $\rightsingularmatk{k}$] =
          \textbf{initializingIncrementalSVD}($\snapshotvec$, $\svdtol$, $k$), i.e., apply
          Algorithm~\ref{al:initializingIncrementalSVD} 
          \RETURN
        \ENDIF
        \STATE $\boldsymbol \projected \assign \leftsingularmatk{k-1}^T\snapshotvec$ 
        \STATE $\projerr \assign \sqrt{\snapshotvec^T\snapshotvec-\projected^T\projected}$
        \STATE $\updatedBasisVector \assign (\snapshotvec - \leftsingularmatk{k-1}\projected)/\projerr$
        \STATE $\middlemat \assign \bmat{\diag{s} & \projected \\ \zero & \projerr}$ 
        \IF{$\projerr < \svdtol$}
          \STATE $\middlemat_{\text{end},\text{end}} \assign 0$
        \ENDIF
        \STATE
        $[\middleleftsingularmat,\middlesingularvalmat,\middlerightsingularmat]
        \assign \text{SVD}(\middlemat)$
        \\ $\#$ SVD update 
        \IF{$\projerr<\svdtol$}
          \STATE $\leftsingularmatk{k} \assign
          \leftsingularmatk{k-1}{\middleleftsingularmat}_{(1:\basisrank_{k-1},1:\basisrank_{k-1})}$,
          \quad
          $\singularvalues_k \assign
          \diag{\middlesingularvalmat_{(1:\basisrank_{k-1},1:\basisrank_{k-1})}}$,\quad
          and $\rightsingularmatk{k} \assign \bmat{\rightsingularmatk{k} &
          \zerobold \\ \zerobold & 1 }
            {\middlerightsingularmat}_{(:,1:\basisrank_{k-1})}$
        \ELSE
          \STATE $\leftsingularmatk{k} \assign \bmat{\leftsingularmatk{k-1} &
          \updatedBasisVector}\middleleftsingularmat$, \quad
          $\singularvaluesk{k} \assign \diag{\middlesingularvalmat}$, \quad and
          $\rightsingularmatk{k}
          \assign \bmat{\rightsingularmatk{k} & \zerobold \\ \zerobold &
          1}\middlerightsingularmat$
        \ENDIF
        \\ $\#$ Neglect small singular values: truncation
        \IF{${\singularvaluesk{k}}_{(\basisrank_{k})} < \svtol$}
          \STATE $\singularvaluesk{k} \assign
          {\singularvaluesk{k}}_{(1:\basisrank_{k}-1)}$, \quad
          $\leftsingularmatk{k} \assign
          {\leftsingularmatk{k}}_{(:,1:\basisrank_{k}-1)}$, \quad
          $\rightsingularmatk{k} \assign 
          {\rightsingularmatk{k}}_{(:,1:\basisrank_{k}-1)}$
        \ENDIF
        \\ $\#$ Orthogonalize if necessary 
        \IF{${\leftsingularmatk{k}}^T_{(:,1)}
        {\leftsingularmatk{k}}_{(:,\text{end})} >
        \min\{\svdtol,\tol\cdot \nspace\}$}
          \STATE $[\Qmat,\Rmat] \assign \qr{\leftsingularmatk{k}}$ 
          \STATE $\leftsingularmatk{k} \assign \Qmat$
        \ENDIF
   \end{algorithmic}
  \end{algorithm}

\subsection{Space--time reduced basis in block
structure}\label{sec:blockstructure}
  Forming the space--time basis in Eq.~\eqref{eq:spacetimebasis} through the the
  Kronecker product in Eq.~\eqref{eq:spacetimebasisvec} requires
  $\nspace\ntime\nreducedspace\nreducedtime$ multiplications. This is
  troublesome, not only because it is computationally costly, but also it
  requires too much memory. Fortunately, a block structure of the space--time
  basis in \eqref{eq:spacetimebasis} is available:
  \begin{equation}\label{eq:blockstructure}
  \boldsymbol{\Phi}_{st} = 
    \bmat{ 
    \spatialbasis\diagonaltemporalbasis{1}{1} & \cdots & \cdots & \cdots &
    \spatialbasis\diagonaltemporalbasis{1}{\nreducedtime} \\
    \vdots & \ddots & \vdots & \rddots & \vdots \\
    \vdots & \ldots & \spatialbasis \diagonaltemporalbasis{k}{\reducedtimeindex} & \ldots &
    \vdots \\
    \vdots & \rddots & \vdots & \ddots & \vdots \\
    \spatialbasis\diagonaltemporalbasis{\ntime}{1} & \cdots & \cdots & \cdots &
    \spatialbasis\diagonaltemporalbasis{\ntime}{\nreducedtime} \\
  } \in \mathbb{R}^{\nspace\ntime\times \nreducedspace\nreducedtime},
  \end{equation}
  where $k$th time step temporal basis matrix,
  $\diagonaltemporalbasis{k}{\reducedtimeindex} \in
  \RR{\nreducedspace\times\nreducedspace}$, is defined as  
  \begin{equation}\label{eq:temporaldiagonal}
    \diagonaltemporalbasis{k}{\reducedtimeindex} \equiv \text{diag}\left (
    \bmat{\temporalbasiselementk{1}{\reducedtimeindex}{k} & \cdots &
    \temporalbasiselementk{\nreducedspace}{\reducedtimeindex}{k} } \right ),
  \end{equation}
  where $\temporalbasiselementk{\reducedspaceindex}{\reducedtimeindex}{k}\in\RR{}$ denotes 
  a $k$th element of $\temporalbasisveck{\reducedspaceindex}{\reducedtimeindex}$. Thanks to this block structure, 
  the space--time reduced order operators, such as $\reducedststatemat$,
  $\reducedstinputvec$, and $\reducedstinitial$ can be formed without explicitly
  forming $\spacetimebasis$. For example, the reduced space--time system matrix,
  $\reducedststatemat$ can be computed, using the block structures, as
  \begin{equation}\label{eq:block-reduced-state-mat}
    \reducedststatemat(\param) = \spacetimebasis^T \ststatemat(\param)
    \spacetimebasis = 
    \bmat{\reducedststatematk{1}{1}(\param) & \cdots & \cdots & \cdots &
    \reducedststatematk{1}{\nreducedtime}(\param) \\
    \vdots & \ddots & \vdots & \rddots & \vdots \\
    \vdots & \cdots & \reducedststatematk{\reducedtimeindex'}{\reducedtimeindex}(\param) & \cdots & \vdots \\
    \vdots & \rddots & \vdots & \ddots & \vdots \\
    \reducedststatematk{\nreducedtime}{1}(\param) & \cdots & \cdots & \cdots &
    \reducedststatematk{\nreducedtime}{\nreducedtime}(\param)
  },
  \end{equation}
  where ($\reducedtimeindex'$,$\reducedtimeindex$)th block matrix,
  $\reducedststatematk{\reducedtimeindex'}{\reducedtimeindex}(\param) \in
  \RR{\nreducedspace\times\nreducedspace}$, $\reducedtimeindex'$,
  $\reducedtimeindex \in \nat{\nreducedtime}$ is defined as
  \begin{equation}\label{eq:ij-block-matrix}
    \reducedststatematk{\reducedtimeindex'}{\reducedtimeindex}(\param) = \sum_{k=1}^{\ntime} \left (
    \diagonaltemporalbasis{k}{\reducedtimeindex'} \diagonaltemporalbasis{k}{\reducedtimeindex} - \timestepk{k}
    \diagonaltemporalbasis{k}{\reducedtimeindex'}
    \reducedstatemat(\param) \diagonaltemporalbasis{k}{\reducedtimeindex}
    \right ) - \sum_{k=1}^{\ntime-1}
    \diagonaltemporalbasis{k+1}{\reducedtimeindex'}
    \diagonaltemporalbasis{k}{\reducedtimeindex}.
  \end{equation}
  Note that the computations of $\diagonaltemporalbasis{k}{\reducedtimeindex'}
  \diagonaltemporalbasis{k}{\reducedtimeindex}$ and
  $\diagonaltemporalbasis{k+1}{\reducedtimeindex'}\diagonaltemporalbasis{k}{\reducedtimeindex}$ are trivial
  because they are diagonal matrix-products whose individual product requires
  $\nreducedspace$ scalar products. Additionally, $\reducedstatemat:\RR{\nparam}
  \rightarrow \RR{\nreducedspace\times\nreducedspace}$, is a reduced order
  system operator that is used for the spatial ROMs, e.g., see
  Eq.~\eqref{eq:galerkin}. This can be pre-computed.  It implies that the
  construction of $\reducedststatemat(\param)$ requires computational cost that
  is a bit larger than the one for the spatial ROM system matrix. The additional
  cost $O(\nreducedspace^2\nreducedtime^2\ntime)$ to the spatial ROM system matrix construction
  is required.

  Similarly, the reduced space--time input vector,
  $\reducedstinputvec\in\RR{\nreducedspace\nreducedtime}$, can be computed as
  \begin{equation}\label{eq:block-reduced-input-vec}
    \reducedstinputvec = \spacetimebasis^T\stinputvec = 
    \bmat{\vdots \\ \reducedstinputveck{\reducedtimeindex} \\ \vdots},
  \end{equation}
  where the $\reducedtimeindex$th block vector, $\reducedstinputveck{\reducedtimeindex}\in\RR{\nreducedspace}$,
  $\reducedtimeindex\in\nat{\nreducedtime}$, is given as
  \begin{equation}\label{eq:j-block-vector}
    \reducedstinputveck{\reducedtimeindex}(\param) = \sum_{k=1}^{\ntime} \timestepk{k}
    \diagonaltemporalbasis{k}{\reducedtimeindex}\reducedinputmat(\param)\inputveck{k}.
  \end{equation}
  Note that $\reducedinputmat: \RR{\nparam} \rightarrow
  \RR{\nreducedspace\times\ninput}$ is used for the spatial ROMs, e.g., see
  Eq.~\eqref{eq:galerkin}. Also, $\reducedinputmat(\param)\inputveck{k}$ needs
  to be computed in the spatial ROMs.
  These can be pre-computed.
  Other operations are
  related to the row-wise scaling with diagonal term,
  $\diagonaltemporalbasis{k}{\reducedtimeindex}$, whose computational cost is
  $O(\nreducedspace\nreducedtime\ntime)$.
  If $\inputveck{k}$ is constant throughout the whole time steps, i.e.,
  $\inputveck{k} = \inputvec$, then
  Eq.~\eqref{eq:j-block-vector} can be further reduced to
  \begin{equation}\label{eq:j-block-vector-constant-input}
    \reducedstinputveck{\reducedtimeindex}(\param) = \left (\sum_{k=1}^{\ntime} \timestepk{k}
    \diagonaltemporalbasis{k}{\reducedtimeindex}\right ) \reducedinputmat(\param)\inputvec,
  \end{equation}
  where you can compute the summation term first, then multiply the diagonal
  term with the precomputed term, $\reducedinputmat\inputvec$, which is not much
  more than the cost for constructing the reduced input vector for the spatial
  ROM. 
  
  Finally, the space--time initial vector,
  $\reducedstinitial\in\RR{\nreducedspace\nreducedtime}$, can be computed as
  \begin{equation}\label{eq:block-reduced-initial-vec}
    \reducedstinitial = \spacetimebasis^T\stinitial = 
    \bmat{\reducedstinitialk{1} \\ \zero \\ \vdots \\ \zero},
  \end{equation}
  where the first block vector, $\reducedstinitialk{1}\in\RR{\nreducedspace}$,
  is given as
  \begin{equation}\label{eq:first-block-initial-vector}
    \reducedstinitialk{1} =
    \left (\sum_{k=1}^{\ntime}\diagonaltemporalbasis{k}{1} \right )
    \reducedinitialstate.
  \end{equation}
  Note that $\reducedinitialstate\in\RR{\nreducedspace}$ is the reduced initial
  condition in the spatial ROM, e.g., see Eq.~\eqref{eq:galerkin}. The
  additional cost to construct the space--time reduced initial vector is
  $O(\nreducedspace)$.
  
  In summary, the block structure in Eq.~\eqref{eq:blockstructure} enables the
  block term expression of the space--time reduced operators, i.e.,
  $\reducedststatemat$, $\reducedstinputvec$, and $\reducedstinitial$, which
  results in a comparable computational cost that is not much more expensive
  than the construction of the spatial reduced operators, i.e.,
  $\reducedstatemat$, $\reducedinputmat$, and $\reducedinitialstate$. This fact
  attracts the desire to use the spatio-temporal ROM rather than spatial ROM
  because the spatio-temporal ROM solving time is much smaller than the
  corresponding spatial ROM.

\section{Error analysis}\label{sec:erroranalysis}
  The error analysis for spatial and spatio-temporal ROMs is presented in this
  section. Section~\ref{sec:errorSpaceROM} presents two error bounds for the
  spatial ROM, while Section~\ref{sec:errorSpacetimeROM} presents error bounds
  for the spatio-termporal ROM.
  
  \subsection{Error analysis for the spatial ROM}\label{sec:errorSpaceROM}
  First, two error bounds will be derived for the spatial ROM
  presented in Section~\ref{sec:spaceROM}.  We define residual function for
  $k$th time step FOM,
  $\resk{k}:\RR{\nspace}\times\RR{\nspace}\rightarrow\RR{\nspace}$, as
  \begin{equation}\label{eq:residual}
    \resk{k}\left(\statek{k},\statek{k-1}\right) \equiv \statek{k} - \statek{k-1} -
    \timestepk{k}\statemat \statek{k} - \timestepk{k}\inputmat\inputveck{k},
  \end{equation}
  which is zero if $\statek{k}$ and $\statek{k-1}$ are FOM solutions from
  Eq.~\eqref{eq:fom}. Here, we drop the parameter dependence for brevity.  Let
  $\approxstatek{k} \in \RR{\nspace}$ be the solution approximation at $k$th
  time step due to the spatial ROM, i.e., $\approxstatek{k} = \refstate +
  \spatialbasis \reducedstatek{k}$. Note that the approximate solutions make
  the following approximate residual function zero:
  \begin{equation}\label{eq:approxres}
    \zero = \approxresk{k}\left(\approxstatek{k},\approxstatek{k-1}\right)
    \equiv
    \approxstatek{k} - \approxstatek{k-1} -
    \timestepk{k}\spatialbasis\spatialbasis^T\statemat \approxstatek{k} -
    \timestepk{k}\spatialbasis\spatialbasis^T\inputmat\inputveck{k}.
  \end{equation}
  Throughout this section, we use the following notations:
  $\diffstatek{k} \equiv \statek{k} - \approxstatek{k}$,
  $\projmat \equiv \spatialbasis\spatialbasis^T$, and
  $\timeintegratoropk{k} \equiv \identityk{\nspace} - \timestepk{k}\statemat$.

  {\theorem\label{theorem:res-basedErrorBoundforSpace} (a residual-based a
  posteriori error bound with the backward Euler time integrator) Let
  $\soperatornormk{k}\in\RR{}$ be a matrix norm of the inverse of the backward
  Euler time integrator operator, i.e., $\soperatornormk{k} \equiv \left \|
  \timeintegratoropk{k}^{-1} \right \|_2$. Then, a residual-based a posteriori
  error bound at $k$th time step is given as
   \begin{equation}\label{eq:res-based-errorbound-for-spatialROM}
     {\displaystyle \left \|\diffstatek{k} \right \|_2 \leq
     \sum_{i=1}^k  \stabilityconstantk{i} \left
     \|\resk{i}\left (\approxstatek{i},\approxstatek{i-1} \right )
     \right \|_2 },
   \end{equation}
   where the stability constants, $\stabilityconstantk{i}\in\RR{}$, are defined
   as $\stabilityconstantk{i} \equiv \prod_{j=i}^{k}\soperatornormk{j}$. 
\vspace*{1em}

  \proof Approximate solutions, $\approxstatek{k}$ and $\approxstatek{k-1}$,
  make the residual nonzero and it can be expanded as 
  \begin{align}\label{eq:res-based-proof}
    - \resk{k} \left(\approxstatek{k}, \approxstatek{k-1}\right) &= \resk{k}
    \left(\statek{k},\statek{k-1}\right) - \resk{k} \left(\approxstatek{k},
    \approxstatek{k-1}\right)  \\ &= \diffstatek{k}  -  \diffstatek{k-1} -
    \timestepk{k} \statemat \diffstatek{k},
  \end{align}
  where we have used the fact that $\resk{k} \left(\statek{k},\statek{k-1}\right) =
  \zero$.  Rearranging terms and inverting the time integrator operator gives
  \begin{equation}\label{eq:rearrange}
    \diffstatek{k} = \timeintegratoropk{k}^{-1} \left [
        \diffstatek{k-1} - \resk{k} 
      \left(\approxstatek{k}, \approxstatek{k-1}\right) \right ].
  \end{equation}
  Taking a norm each side, applying triangle inequality and H{\" o}lder's
  inequality, we obtain the following one-time step bound:
  \begin{equation}\label{eq:one-timestep-bound}
    \left \| \diffstatek{k} \right \|_2 \leq \soperatornormk{k} \left ( \left \|
    \diffstatek{k-1} \right \|_2 + \left \| \resk{k}\left(\approxstatek{k},
    \approxstatek{k-1}\right) \right \|_2 \right )
  \end{equation}
  Assuming $\approxstatek{0} = \statek{0}$, which can be achieved by setting
  $\refstate = \statek{0}$, and applying \eqref{eq:one-timestep-bound}
  recursively, we get the claimed bound, i.e.,
  \eqref{eq:res-based-errorbound-for-spatialROM}.
    \hspace*{1em}\endproof} 
\vspace*{1em}

  It is easy to see that the error bound in
  \eqref{eq:res-based-errorbound-for-spatialROM} is exponentially increasing
  with respect to time because of the summation and product appeared in the
  definition of stability constants, i.e., $\stabilityconstantk{i}$.  For
  $\ell_2$ induced matrix norm, $\soperatornormk{k}$ is the reciprocal of the
  smallest singular value of $\timeintegratoropk{k}$. Also, the error bound in
  \eqref{eq:res-based-errorbound-for-spatialROM} allows any approximate
  solution, that is, $\approxstatek{k}$ does not need to come from the spatial
  ROM solution. The next theorem, however, shows an error bound for a specific
  case, i.e., the error bound for the spatial ROM solutions. 

  {\theorem\label{theorem:specificErrorBoundforSpace} (a spatial ROM-specific a
  posteriori error bound with the backward Euler time integrator) 
  Let $\serrorboundcoefk{k}\in\RR{}$ be defined as $\serrorboundcoefk{k} \equiv
  1/\left (1-\timestepk{k}\|\statemat \|_2\right )$. Also, assume that the
  timestep, $\timestepk{k}$, is sufficiently small, i.e., $\timestepk{k} <
  1/\|\statemat \|_2$. Then, a spatial ROM-specific a posteriori error bound at
  $k$th time step is given as
   \begin{equation}\label{eq:specific-errorbound-for-spatialROM}
     {\displaystyle \left \|\diffstatek{k}\right \|_2 \leq 
      \sum_{i=1}^k  \sconstantserrork{i} \left \|
     \remainveck{i} \right \|_2 },
   \end{equation}
  where the stability constants, $\sconstantserrork{i}\in\RR{}$, are defined as 
  $\sconstantserrork{i} \equiv \prod_{j=i}^{k}\serrorboundcoefk{j}$, 
  and $\remainveck{k}\in\RR{\nspace}$ be defined as $\remainveck{i} \equiv
  \timestepk{i}\left (\statemat\approxstatek{i} +
  \inputmat\inputveck{i}\right)$.  

  \proof Substracting Eq.~\eqref{eq:residual} by Eq.~\eqref{eq:approxres} gives  
  \begin{align}\label{eq:specific-error-proof}
    \zero &= \resk{k} \left(\statek{k}, \statek{k-1}\right) - \approxresk{k}
             \left(\approxstatek{k}, \approxstatek{k-1}\right) \\
    &= \diffstatek{k} - \diffstatek{k-1} - \timestepk{k} \statemat \statek{k} +
    \timestepk{k} \projmat\statemat \approxstatek{k} - \timestepk{k} 
    \projerrormat \inputmat\inputveck{k},
  \end{align}
  where the projection error matrix,
  $\projerrormat\in\RR{\nspace\times\nspace}$, is defined as $\projerrormat
  \equiv
  \identityk{\nspace} - \projmat$.  Rearranging terms, adding and subtracting
  $\timestepk{k}\statemat\approxstatek{k}$, and taking a norm with triangle
  inequality and H{\"o}lders's inequality give
  \begin{equation}\label{eq:specific-rearrange}
    \left \|\diffstatek{k} \right \|_2 \leq  \left \| \diffstatek{k-1} \right \|_2
    + \timestepk{k} \left \| \statemat \right \|_2 \left \| \diffstatek{k} \right \|_2
    + \|\projerrormat\|_2 \left \|\remainveck{k} \right \|_2
  \end{equation}
  Rearranging terms again and dividing by $1-\timestepk{k}\|\statemat\|_2$,
  using the assumption of $\timestepk{k} < 1/\|\statemat \|_2$, give
  \begin{equation}\label{eq:specific-one-timestep-bound}
    \left \| \diffstatek{k} \right \|_2 \leq
    \serrorboundcoefk{k} \left ( \left \| \diffstatek{k-1} \right \|_2 +  \|
    \remainveck{k} \|_2 \right ),
  \end{equation}
  where we used the fact that $\|\projerrormat\|_2 = 1$.
  Assuming $\approxstatek{0} = \statek{0}$, which can be achieved by setting
  $\refstate = \statek{0}$, 
  and applying \eqref{eq:specific-one-timestep-bound}
  recursively, the claimed bound is obtained, i.e.,
  \eqref{eq:specific-errorbound-for-spatialROM}.
    \hspace*{1em}\endproof} 
\vspace*{1em}

  As in the residual-based error bound in
  Theorem~\ref{theorem:res-basedErrorBoundforSpace}, the error bound in
  Theorem~\ref{theorem:specificErrorBoundforSpace} is also exponentially
  increasing with respect to time because of the summation and product appeared
  in the definition of stability constants, i.e., $\sconstantserrork{i}$.
  However, it is easier to see from
  \eqref{eq:specific-errorbound-for-spatialROM} that the effect of exponential
  growth is degraded as the time step decreases. It is because
  $\sconstantserrork{i}$ becomes closer to one as the time step decreases.

  \subsection{Error analysis for the spatio-temporal
  ROM}\label{sec:errorSpacetimeROM}
  Now, we turn our attention to the error bound for the space--time ROM
  solutions. A residual-based error bound will be derived. For the
  backward Euler time integrator, the space--time residual function,
  $\stres:\RR{\nspace\ntime} \rightarrow \RR{\nspace\ntime}$, is
  defined as
  \begin{equation}\label{eq:stresidual}
    \stres(\ststate) \equiv \ststatemat\ststate - \stinputvec - \stinitial, 
  \end{equation}
  where the parameter dependence is dropped. We define the space--time infinity
  norm, $\|\cdot\|_\infty: \RR{\nspace\ntime} \rightarrow \RR{}$, as
  \begin{equation}\label{eq:spacetimeinfnorm}
    \left \| \ststate \right \|_\infty \equiv \max_{k\in\nat{\ntime}} \left \|
    \statek{k} \right \|_2
  \end{equation}
  Throughout this section, the following notation is used:
  $\diffststate \equiv \ststate - \approxststate$.

  {\theorem\label{theorem:res-basedErrorBoundforSpaceTime} (a space--time
  residual-based a posteriori error bound with the backward Euler time
  integrator) 
  The space--time residual-based
  a posteriori error bound is given as
   \begin{equation}\label{eq:res-based-errorbound-for-spacetimeROM}
     {\displaystyle \left \| \diffststate \right \|_\infty \leq \sqrt{\ntime}
     \left \| \left ( \ststatemat \right )^{-1} \right \|_2
     \left \| \stres\left (\approxststate \right ) \right \|_\infty }.
   \end{equation}

  \proof Approximate space--time solution, $\approxststate$, 
  makes the space--time residual nonzero and it can be expanded as 
  \begin{align}\label{eq:res-based-proof}
    \stres \left( \approxststate \right) &= \stres \left( \approxststate \right) 
    - \stres \left( \ststate \right) \\
    &= - \ststatemat \diffststate,
  \end{align}
  where we have used the fact that $\stres \left(\ststate\right) = \zero$.
  Inverting the space--time operator, taking $\ell_2$ norm and H{\"o}lders'
  inequality, and squaring both sides gives
  \begin{equation}\label{eq:stres-rearrange}
    \left \| \diffststate \right \|_2^2 \leq \left \| \left (\ststatemat
    \right)^{-1} \right \|_2^2 \left \| \stres \left ( \approxststate \right )
    \right \|_2^2.
  \end{equation}
  Note that Inequality~\eqref{eq:stres-rearrange} can be re-written as
  \begin{equation}\label{eq:rewritten-inequality}
    \sum_{k=1}^{\ntime} \left \| \diffstatek{k} \right \|_2^2 
    \leq \left \| \left (\ststatemat \right)^{-1} \right \|_2^2 
    \sum_{k=1}^{\ntime} \left \| \resk{k} \left ( \approxstatek{k},
    \approxstatek{k-1}  \right )  \right \|_2^2.
  \end{equation}
  Due to the norm equivalence relations, i.e.,  $\| \dummyvec \|_\infty \leq \|
  \dummyvec \|_2 \leq \sqrt{N} \| \dummyvec \|_\infty$ for a vector
  $\dummyvec\in\RR{N}$, we have
  \begin{equation}\label{eq:infinity-bound}
    \max_{k\in\nat{\ntime}}\left \| \diffstatek{k} \right \|_2 \leq
    \sqrt{\ntime} \left \| \left (\ststatemat \right)^{-1} \right \|_2
    \max_{k\in\nat{\ntime}} \left \| \resk{k} \right \|_2,
  \end{equation}
  which is equivalent to the claimed bound in
  \eqref{eq:res-based-errorbound-for-spacetimeROM}.
    \hspace*{1em}\endproof} 
\vspace*{1em}

\section{The neutron transport equation}\label{sec:NeutronTransport}

\subsection{Boltzmann transport equation}\label{sec:transport}
The Boltzmann equation for the neutron flux function,
$\psi(\position,\energy,\direction,\timevar):
\RR{3}\times\RR{}\times\RR{2}\times\RR{} \rightarrow \RR{}$, is formulated as
\begin{equation}\label{xport0}
  \frac{1}{\speed(\energy)} \frac{\partial \psi}{\partial \timevar} +
  \direction\cdot\nabla\psi + \crosssection(\position,\energy)\psi =
  \int_0^{\infty} \int_{4\pi}
  \scatter(\position,\energyin\rightarrow\energy,\directionin\cdot\direction)
  \psi(\position,\energyin,\directionin,\timevar)d\directionin d\energyin +
  \source(\position,\energy,\direction,\timevar),
\end{equation}
where $\position\in\RR{3}$ denotes a position vector, $\energy\in\RR{}$ denotes
energy, and $\direction\in{\cal{S}}^2$ (the unit sphere in $\RR{2}$) denotes a
directional vector.  The speed of the neutron is a function of energy, i.e.,
$\speed:\RR{}\rightarrow\RR{}$. The cross-sectional area of a target nucleus is
a function of position and energy, i.e.,
$\crosssection(\position,\energy):\RR{3}\times\RR{}\rightarrow\RR{}$.  The
scattering cross-sectional area is denoted as $\scatter$ and an external source
function is denoted as $\source(\position,\energy,\direction,\timevar):
\RR{3}\times\RR{}\times\RR{2}\times\RR{}\rightarrow\RR{}$.  The spatial domain
is the box ${\cal D} \equiv \{ r=(x,y,z) | a_x \leq x \leq b_x, a_y \leq y \leq
b_y, \mbox{ and } a_z \leq z \leq b_z\}$, and the spatial gradient is denoted as
$\nabla \psi \equiv (\partial \psi/\partial x, \partial \psi/\partial y,
\partial \psi/\partial z)$. We also assume that $\int_{{\cal S}^2} d\direction =
1$ as in Lewis and Miller \cite{LM93}.

Boundary conditions must also be specified to make  (\ref{xport0})
well-posed.  Various options include a reflecting 
condition on a face, or a Dirichlet condition in which the incident flux
is specified on a face.  For simplicity, we will consider only the
latter case.  Namely, we will consider vacuum boundary conditions of the form
\begin{equation}
\psi(\position,\direction,\energy, t) = 0 \mbox{ for all $\position \in
  \partial {\cal D}$ and $\direction \in {\cal S}^2$ with }
\vec{n}(\position) \cdot \direction < 0 ,
\label{dirbcs}
\end{equation}
where $\vec{n}(\position)$ is the outward pointing unit normal at
$\position \in \partial{\cal D}$.

A semi-discretization of (\ref{xport0}) can be obtained
using a {\it multigroup\/} discretization of the energy $\energy$ (see, {\em
e.g.}, \cite{LM93}).  In the multigroup 
approach, the energy $\energy$ is restricted to a finite interval
partitioned into subintervals, or ``groups'':
\begin{equation}
\energy_{max} = \energy_0 > \energy_1 > \cdots > \energy_G = \energy_{min} . 
\end{equation}
The equation (\ref{xport0}) is then averaged over each group $\energy_g < \energy
< \energy_{g-1}$ and the cross-sections $\sigma$ and $\sigma_s$ are
approximated by certain ``flux-weighted averages'' to maintain
linearity.  
This yields the following semi-discretization of
(\ref{xport0}):
\begin{eqnarray}
& & \frac{1}{\speed_g} \frac{\partial
    \psi_g(\position,\direction,\timevar)}{\partial \timevar} + \direction
  \cdot \nabla \psi_{g} (\position,\direction,\timevar) + \sigma_g(\position)
  \psi_g(\position,\direction,\timevar) = \label{scalar_semidiscrete} \\ &
  & \hspace{5em} \sum_{g'=1}^G \int_{{\cal S}^2}
  \sigma_{s,g,g'}(\position,\direction' \cdot \direction)
  \psi_{g'}(\position,\direction',\timevar) d\direction' +
  q_{g}(\position,\direction,\timevar), \nonumber
\end{eqnarray}
for $g=1,\cdots,G$, where
$\psi_g(\position,\direction) \equiv \int_g \psi(\position,\direction,\energy) d\energy$ and
$q_g(\position,\direction) \equiv \int_g q(\position,\direction,\energy) d\energy$, with $\int_g d\energy =
\int_{\energy_g}^{\energy_{g-1}} d\energy$.  

When solving (\ref{scalar_semidiscrete}), for each $g$ the flux
$\psi_g(\position,\direction)$ is expanded in surface harmonics according to 
\[
\psi_g(\position,\direction) = \sum_{n=0}^{\infty} \sum_{m=-n}^{n}
 \phi_{g,n,m}(\position) \; Y_n^m(\direction) .
\]
Here, $Y_n^m(\direction)$ is a surface harmonic defined by 
\[
Y_n^m(\direction) = a_n^m P_n^{|m|}(\xi) \tau_m(\varphi) ,
\]
where $\direction=(\mu,\eta,\xi) = (\sin \theta \cos \varphi, \sin \theta \sin
\varphi, \cos \theta)$, $P_n^{|m|}$ is an {\em associated Legendre
polynomial} \cite{Liboff80}, and 
\[
\tau_m(\varphi) = \left\{ \begin{array}{l}
         \cos m\varphi \mbox{, if $m \geq 0$, and} \\
         \sin |m|\varphi \mbox{, if $m < 0$} .
         \end{array} \right.
\]
The constants $a_n^m$ are defined by
\[
a_n^m = \left[
\frac{2(2n+1)(n-|m|)!}{(1+\delta_{m0}) (n+|m|)!}\right]^{1/2},
\]
where $\delta_{n,n'}$ is the {\em Kronecker delta}, and
\[
  \phi_{g,n,m}(\position)
  \equiv \int_{{\cal S}^2} \psi_{g}(\position,\direction) Y_n^m(\direction)d\direction ,
\]
is the $(n,m)^{th}$ moment of $\psi$.
We have
\[
\int_{{\cal S}^2} Y_n^m(\direction) \; Y_{n'}^{m'}\!(\direction) \; d\direction
= \delta_{n,n'} \; \delta_{m,m'}
\]
for all $n,n'=0,1,\cdots$, and $|m|\leq n, |m'| \leq n'$.
The source $q_g$ is similarly expanded. 

Given $\psi_g$ in the above form, one is able to rewrite the scattering
integral in the form
\begin{equation}
\int_{{\cal S}^2} \sigma_{s,g,g'}(\position,\direction' \cdot \direction)
\psi_{g'}(\position,\direction')d\direction' 
=\sum_{n=0}^{\infty} \sigma_{s,g,g',n}(\position) \sum_{m=-n}^{n}
\phi_{g',n,m}(\position) Y_n^m(\direction), 
\label{scattering}
\end{equation}
where the $\sigma_{s,g,g',n}$ are given by
\[
\sigma_{s,g,g',n}(\position) \equiv 
        \frac{1}{2} \int_{-1}^{1} \sigma_{s,g,g'}(\position,\mu_0)
        P_n(\mu_0) d\mu_0 , 
\]
and where $\mu_0$ is the cosine of the scattering angle.
The infinite series in (\ref{scattering}) is truncated to a finite
number of terms, with a maximum value $N_s$ for $n$.  Thus, we can
write the multigroup equations as
\begin{eqnarray}\label{eq:multigroupequation}
& &
\frac{1}{\speed_g} \frac{\partial
  \psi_g(\position,\direction)}{\partial \timevar} + 
\direction \cdot \nabla \psi_{g} (\position,\direction) + \sigma_g(\position) 
\psi_g(\position,\direction) =  \label{multigrp_eqns} \\
& & \hspace{5em} \sum_{g'=1}^G 
\sum_{n=0}^{N_s} \sigma_{s,g,g',n}(\position) \sum_{m=-n}^{n}
\phi_{g',n,m}(\position) Y_n^m(\direction) + q_{g}(\position,\direction), \nonumber 
\end{eqnarray}
for $g=1,\cdots,G$.

\subsection{Spatial and directional discretization of the 3-D Problem}

In previous work \cite{BiBr2009}, we derived a matrix version of the
well-known {\em simple corner balance} (\SCB) discretization scheme for the
1-D slab problem analogous to (\ref{xport0})-(\ref{dirbcs}). See
\cite{Adam97} for more details about the {\SCB} method.
This matrix formalism can easily be extended to 3-D problems, and we give
a brief overview here.

The angular variable $\direction$ is discretized using a quadrature
rule.  The specific quadrature rules we consider for approximating
integrals on ${\cal S}^2$ employ the standard symmetry
assumptions. Following Carlson and Lathrop \cite{CaLa68}, we consider
quadrature rules of the form
\begin{equation}
\int_{{\cal S}^2} \psi(\direction)d\direction \approx \sum_{\ell=1}^{L}
w_{\ell} \psi(\direction_{\ell}),
\end{equation}
where $\direction_{\ell} \equiv (\mu_{\ell}, \eta_{\ell}, \xi_{\ell})$,
for all $\ell = 1,\ldots,L$, with $L = \nu(\nu+2)$ and $\nu$ is the
number of direction cosines ($\nu = 2,4,6,\ldots$).

For the spatial discretization, we begin by considering the
mono-energetic steady-state Boltzmann equation
\begin{equation}
\left\{ \begin{array}{l}
\direction \cdot \nabla \psi + \sigma \psi = f \mbox{ in } {\cal D}, \\
\psi(\position) = 0 \mbox{ for all $\position \in \partial{\cal D}$ with }
\vec{n}(\position) \cdot \direction < 0, \end{array} \right. \label{hypera}
\end{equation}
where $\direction=(\mu,\eta,\xi) \in {\cal S}^2$ is fixed and equal to
one of the above quadrature points (although we suppress the $\ell$
subscript to simplify notation), ${\cal D}$ is the spatial domain
defined earlier, and $\vec{n}(\position)$ is the outward pointing unit
normal at $\position \in \partial{\cal D}$.  The functions $f$ and
$\sigma$ are assumed known. The spatial domain ${\cal D}$ is
discretized into zones in the natural way, defining
\begin{eqnarray*}
& & \Delta x_i = x_i - x_{i-1} \mbox{ for } i=1,\ldots,M, \;\;
    \Delta y_j = y_j - y_{j-1} \mbox{ for $j=1,\ldots,J$, and} \\
& & \Delta z_k = z_k - z_{k-1} \mbox{ for } k=1,\ldots,K, 
\end{eqnarray*}
and define $\position_{ijk} = (x_i,y_j,z_k)$.  Also define $\Delta
\position_{ijk} \equiv \Delta x_i \Delta y_j \Delta z_k$.  We will
view the {\SCB} method here as a zone-centered discretization without
the use of $8 \times 8$ superzones in 3D as is normally done with
{\SCB}. Thus, the parameters $M$, $J$, and $K$ must all be even
numbers.  Assume that $\sigma$ and $f$ have constant values on each
{\em zone\/}
\[
{\cal Z}_{ijk} \equiv \{\position | x_{i-1} < x < x_i, y_{j-1} < y < y_j,
                          z_{k-1} < z < z_k\}, 
\]
denoted by $\sigma_{ijk}$ and $f_{ijk}$, respectively. We use
$\psi_{ijk}$ to denote the approximation to $\psi(\position_{ijk})$,
the true solution at $\position_{ijk}$.  Following the development
given in \cite{BiBr2009}, there are $MJK$ unknowns
$\psi_{ijk}$, and  $MJK$ equations.

Writing the discretized system in matrix notation, we first have the
discrete flux vector and right hand side
\[
  \boldsymbol{\Psi} \mbox{ and } \boldsymbol{F} \in \RR{MJK},
\]
defined for all zones ordered by $i$ first, then $j$, and finally
$k$. Next, define the diagonal matrix 
\begin{equation}
  \boldsymbol{\Sigma} \equiv \mbox{diag}(\sigma_{111},\cdots,\sigma_{MJK}).
\label{Sigt}
\end{equation}
The SCB discretization of the $\direction \cdot \nabla$ operator then
results in matrices $\boldsymbol{C}_x$, $\boldsymbol{C}_y$, and $\boldsymbol{C}_z \in \RR{MJK \times MJK}$,
similar in form to the $\boldsymbol{G}_j$ matrices in \cite{BiBr2009}, but permuted
because of the cell-centered ordering, and has the form
\begin{equation}
  \boldsymbol{C} \equiv \mu \boldsymbol{C}_x + \eta \boldsymbol{C}_y + \xi
  \boldsymbol{C}_z \approx \direction \cdot \nabla.
\label{Ctimes}
\end{equation}
While not explicitly noted, the $\boldsymbol{C}_x$, $\boldsymbol{C}_y$, and
$\boldsymbol{C}_z$ matrices in
this approximation also depend on the particular octant of ${\cal
  S}^2$ the variable $\direction$ is in.  Putting (\ref{Sigt}) and
(\ref{Ctimes}) together we have (and adding the quadrature point and
group dependence) the matrix representation of the discrete version of
(\ref{hypera}) can be written as
\begin{equation}
  \boldsymbol{H}_{g,\ell} \boldsymbol{\Psi}_{g,\ell} = \boldsymbol{F}_{g,\ell} 
\mbox{, where }
  \boldsymbol{H}_{g,\ell} \equiv \boldsymbol{C}_{\ell} + \boldsymbol{\Sigma}_g .
\label{HellDdef}
\end{equation}

\subsection{The discrete ordinates method}

Continuing the matrix development of the overall
discretization of (\ref{xport0}), we begin by defining discretized
representations of the operations of taking moments of the flux.
As operators on zone-centered vectors, these are easily seen to be
given by the $MJK \times LMJK$ matrices  
\begin{equation}
  \boldsymbol{L}_{n,m} \equiv (\boldsymbol{l}_{n,m} \boldsymbol{W}) \tensor \boldsymbol{I}_{MJK}  \label{Lnmdef}
\end{equation}
where
\[
  \boldsymbol{l}_{n,m} \equiv ( Y_n^m(\direction_1), Y_n^m(\direction_2), \cdots ,
           Y_n^m(\direction_{L}) )  \label{lnmdef}
\mbox{ and } \boldsymbol{W} \equiv \mbox{diag}(w_1,\cdots,w_L).
\]
If the vector $\boldsymbol{\Psi}_g$ approximates $\psi_g(\position,\direction)$,
then $\boldsymbol{L}_{n,m}
\boldsymbol{\Psi}_g$ will approximate the $(n,m)^{\rm th}$ moment of
$\psi_g(\position,\direction)$, namely $\phi_{g,n,m}(\position)$. Similarly, we define the
$LMJK \times MJK$ matrices
\begin{equation}
  \boldsymbol{L}_{n,m}^{+} \equiv \boldsymbol{l}_{n,m}^T \tensor \boldsymbol{I}_{MJK}. \label{Lnm+def}
\end{equation}
For a vector $\boldsymbol{\Phi}$ approximating $\phi(\position)$,
$\boldsymbol{L}_{n,m}^{+}
\boldsymbol{\Phi}$ will
approximate $Y_n^m(\direction) \phi(\position)$. We also will find it useful to
define the grouped matrices $\boldsymbol{L}_n$ and $\boldsymbol{L}_n^+$, where
\[
  \boldsymbol{L}_n = \left( \begin{array}{c} \boldsymbol{L}_{n,-n} \\ \vdots \\
    \boldsymbol{L}_{n,n}
                    \end{array} \right)
\mbox{ and }
\boldsymbol{L}^+_n = \left(\boldsymbol{L}^+_{n,-n},\cdots,\boldsymbol{L}^+_{n,n} \right) ,
\]
and also the further grouped matrices
\[
  \boldsymbol{L}^N = \left( \begin{array}{c} \boldsymbol{L}_{0} \\ \vdots \\
    \boldsymbol{L}_{N}
                    \end{array} \right)
\mbox{ and }
\boldsymbol{L}^{N,+} = \left(\boldsymbol{L}^+_{0},\cdots,\boldsymbol{L}^+_{N} \right) ,
\]
Given an $N = N_s$, the number of terms in the scattering kernel, we
will assume that the quadrature rule is symmetric through the origin
(see remarks above) and such that the spherical harmonics of order
$N_s$ and less satisfy
\[
\sum_{\ell = 1}^{L} Y_n^m(\direction_{\ell}) Y_{n'}^{m'}(\direction_{\ell}) =
\delta_{n,n'} \delta_{m,m'}
\]
for all $0 \leq n,n' \leq N_s$, $|m|\leq n$, and $|m'| \leq n'$. This
can be written more compactly as 
\begin{equation}
  \boldsymbol{L}^{N_s} \boldsymbol{L}^{N_s,+} = \boldsymbol{I}_{(N_s+1)^2}
  \tensor \boldsymbol{I}_{MJK}.
\label{L_Lp_eqn}
\end{equation}
To represent the source term, define the zone-centered vector $\boldsymbol{Q}
\equiv (q_{ijk\ell}) \in \RR{LMJK}$, where $q_{ijk\ell} \equiv
q(r_{ijk},\direction_{\ell})$. Next, let
\begin{eqnarray}
  \boldsymbol{\Sigma}_{s,g,g',n} &\equiv& \boldsymbol{I}_{2n+1} \tensor
  \hat{\boldsymbol{\Sigma}}_{s,g,g',n},
          \hspace{1em} \mbox{where } \label{gam1hat} \\
           & & \hat{\boldsymbol{\Sigma}}_{s,g,g',n} \equiv {\rm
diag}(\sigma_{s,g,g',n,111},\ldots, \sigma_{s,g,g',n,MJK}) , \hspace{1em}
           n=0,1,\ldots \mbox{, and } \nonumber \\ 
  \bar{\boldsymbol{\Sigma}} &\equiv& \boldsymbol{I}_{L} \tensor
  \boldsymbol{\Sigma}. \nonumber
\end{eqnarray}
Using the above matrices, define the matrix $\boldsymbol{H}_g$ by
\begin{equation}
  \boldsymbol{H}_g \equiv \mbox{diag}(H_{g,1},\ldots,H_{g,L}).
\label{Hdef}
\end{equation}
If we assume only $N_s + 1$ terms in the scattering operator, then the
complete discretization of (\ref{xport0})--(\ref{dirbcs}) can be
written in the compact form
\begin{equation}
  v_g^{-1} \dot{\boldsymbol{\Psi}}_g +  \boldsymbol{H}_g \boldsymbol{\Psi}_g =
  \bar{\boldsymbol{Z}} \sum_{g'=1}^{G} \sum_{n=0}^{N_s} \boldsymbol{L}_n^{+}
         \boldsymbol{\Sigma}_{s,g,g',n} \boldsymbol{L}_n \bar{\boldsymbol{S}}
         \boldsymbol{\Psi}_{g'} 
         + \boldsymbol{Q}_g,  \; \; g=1,\cdots,G .
\label{matddeqns}
\end{equation}
Next, if we define
\[
  {\bf \Psi} \equiv \left( \begin{array}{c} \boldsymbol{\Psi}_1 \\
    \boldsymbol{\Psi}_2 \\
                                          \vdots \\
  \boldsymbol{\Psi}_G \end{array} \right),
\; \;
{\bf \Psi_B} \equiv \left( \begin{array}{c} \boldsymbol{\Psi}_{B,1} \\
  \boldsymbol{\Psi}_{B,2} \\
                                            \vdots \\
\boldsymbol{\Psi}_{B,G} \end{array} \right),
\; \;
{\bf F} \equiv \left( \begin{array}{c} \boldsymbol{F}_1 \\
  \boldsymbol{F}_2 \\
                                       \vdots \\
\boldsymbol{F}_G \end{array} \right),
\; \;
{\bf Q} \equiv \left( \begin{array}{c} \boldsymbol{Q}_1 \\
  \boldsymbol{Q}_2 \\
                                       \vdots \\
\boldsymbol{Q}_G \end{array} \right),
\; \;
{\bf \Sigma}_s \equiv \left( \begin{array}{ccc}
  \boldsymbol{\Sigma}_{s,11}^{N_s} & \cdots & \boldsymbol{\Sigma}_{s,1G}^{N_s} \\
               \vdots              & \ddots & \vdots              \\
               \boldsymbol{\Sigma}_{s,G1}^{N_s} & \cdots & \boldsymbol{\Sigma}_{s,GG}^{N_s}
                \end{array} \right) ,
\]
where $\boldsymbol{\Sigma}_{s,gg'}^{N_s} \equiv \mbox{diag} \left(
\boldsymbol{\Sigma}_{s,g,g',0},\cdots, \boldsymbol{\Sigma}_{s,g,g',N_s} \right)$, and define
${\bf H} \equiv
\mbox{diag}(\boldsymbol{H}_1,\boldsymbol{H}_2,\cdots,\boldsymbol{H}_G)$, ${\bf L}^{+} \equiv
\boldsymbol{I}_G \tensor \boldsymbol{L}^{N_s,+}$, ${\bf L} \equiv
\boldsymbol{I}_G \tensor \boldsymbol{L}^{N_s}$, ${\bf V}
\equiv \boldsymbol{V} \tensor \boldsymbol{I}_{LMJK}$, with $\boldsymbol{V} \equiv
\mbox{diag}(v_1,\cdots,v_G)$, then (\ref{matddeqns}) can be written as
\begin{equation}
{\bf V}^{-1} \dot{\bf \Psi} + {\bf H} {\bf \Psi} = {\bf L}^{+} {\bf
  \Sigma}_s {\bf L} {\bf \Psi} + {\bf Q}, \label{DAE}.
\end{equation}
Finally, writing \eqref{DAE} in the form of systems
(Eqs.~\eqref{eq:Dynamicalstate} and \eqref{eq:Dynamicaloutput}) gives
\begin{eqnarray*}
\dot{\bf \Psi}(t) &=& {\bf V} \left( {\bf L}^{+} {\bf
  \Sigma}_s {\bf L} - {\bf H}  \right) {\bf \Psi} +
    {\bf V} {\bf Q}(t), \quad {\bf \Psi}(0) = {\bf \Psi}_0 , \\
    \boldsymbol{R}(t) &=& {\bf D}^T {\bf \Psi}(t),
\end{eqnarray*}
with $\boldsymbol{R}(t) \in \RR{}$ representing the response of the flux ${\bf \Psi}(t)$
integrated over a region in phase space and the matrix $\bf D$
performs the integration.

\section{Numerical results}\label{sec:numericalresults}
\begin{figure}[h]
\centering
  \subfloat[Geometry for Example 1 in Section~\ref{sec:ex1}\label{fig:ex1}]
  {\includegraphics[width=0.4\textwidth]{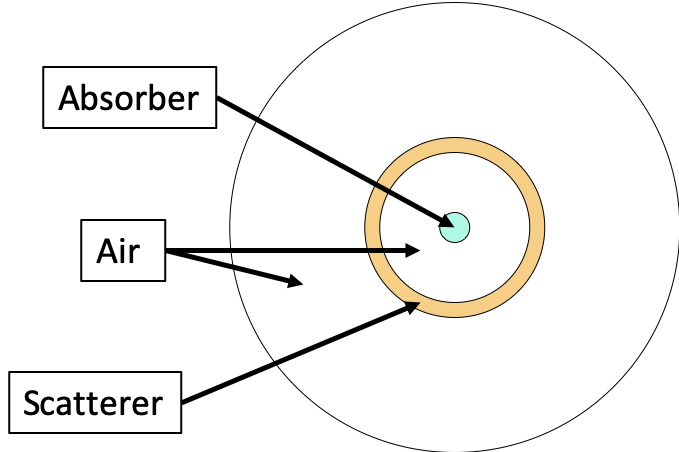}}\hfill
  \subfloat[Geometry for Example 2 in Section~\ref{sec:ex2}\label{fig:ex2}] 
  {\includegraphics[width=0.4\textwidth]{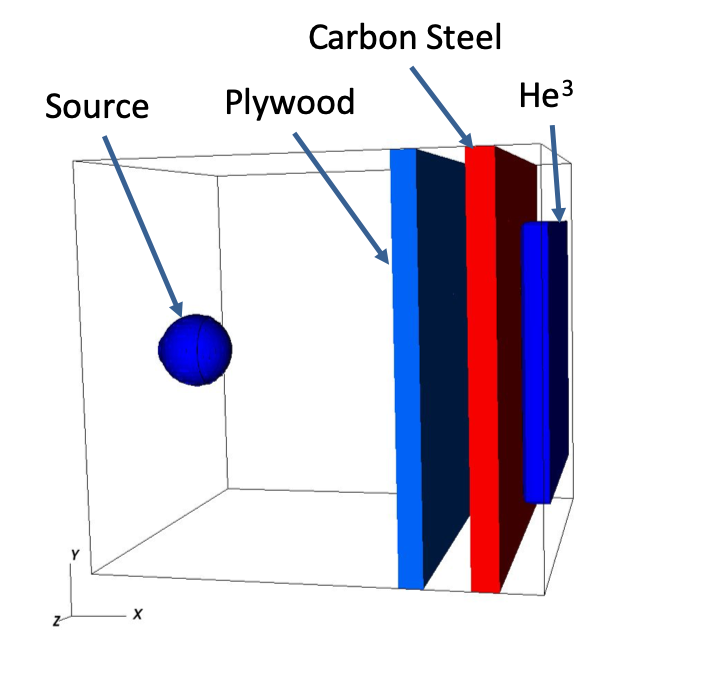}}
\caption{Geometry description of numerical examples} 
\label{fig:examples}
\end{figure}

  We present the performance of the space--time ROM applied to the Boltzmann
  particle transport equation.  We consider two 3D neutron particle simulation
  examples with two different geometries, see Figure~\ref{fig:examples}.  The
  full order model simulations are done by
  ARDRA\footnote{https://computing.llnl.gov/projects/ardra-scaling-up-sweep-transport-algorithms},
  the LLNL production code for the transport sweep algorithms
  \cite{bihari2009linear}.  The space--time ROM is implemented within ARDRA
  source code and the reduced bases are generated by
  libROM\footnote{https://computing.llnl.gov/projects/librom-pod-based-reduced-order-modeling},
  i.e., the LLNL reduced order basis generation codes \cite{osti_1505575}.  The libROM
  can be obtained from the following github page:
  https://github.com/LLNL/libROM.  The Boltzmann particle transport equation and
  its numerical discretization are described in
  Section~\ref{sec:NeutronTransport}.  Each example description is detailed in
  Sections~\ref{sec:ex1} and \ref{sec:ex2}. 

  All the simulations in this numerical section use RZTopaz in Livermore
  Computing Center\footnote{https://hpc.llnl.gov/hardware/platforms/RZTopaz}, on
  Intel Xeon CPUs with 128 GB memory, peak TFLOPS of 928.9, and peak single CPU
  memory bandwidth of 77 GB/s.

\subsection{Example 1: a symmetric case}\label{sec:ex1}
  The first example solves the Boltzmann particle transport equation on a 3D
  cartesian mesh.  The mesh is $20\times20\times20$, resulting in $8,000$
  spatial zones.  There are eighty angular directions and seventeen energy
  groups. An absorber is located at the center and the second shell is scatterer
  as described in Figure~\ref{fig:ex1}. The neutron source is $14.1$ $MeV$,
  which is in the 2nd energy group. The source is constant and the final
  simulation time step is at $\finaltime=40$ $nsec$ with a uniform time step
  $\timestep = 1$ $nsec$.  As a result, there are $10,880,000$ degrees of
  freedom in space and $435,200,000$ degrees of freedom in space--time.
  The full order model simulation uses $8$ cores in RZTopaz and takes 22.5
  seconds, resulting in the CPU time of around 3 minutes.

\begin{figure}[h!]
  \centering
    \subfloat[Neutron flux at $t=0$ $nsec$\label{fig:initial1}]
    {\includegraphics[width=0.27\textwidth]{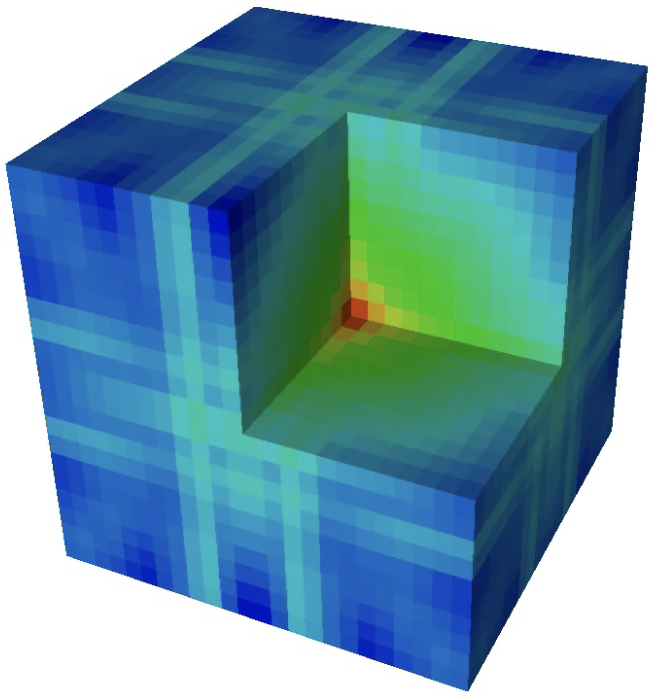}}\hfill
    \subfloat[Neutron flux at $t=40$ $nsec$\label{fig:final1}]
    {\includegraphics[width=0.27\textwidth]{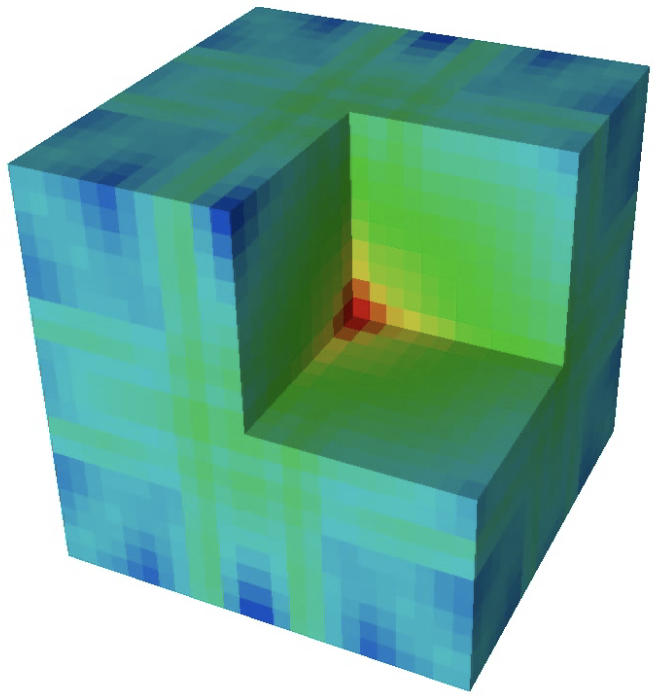}}\hfill
    \subfloat[relative error  \label{fig:relerror1}] 
    {\includegraphics[width=0.4\textwidth]{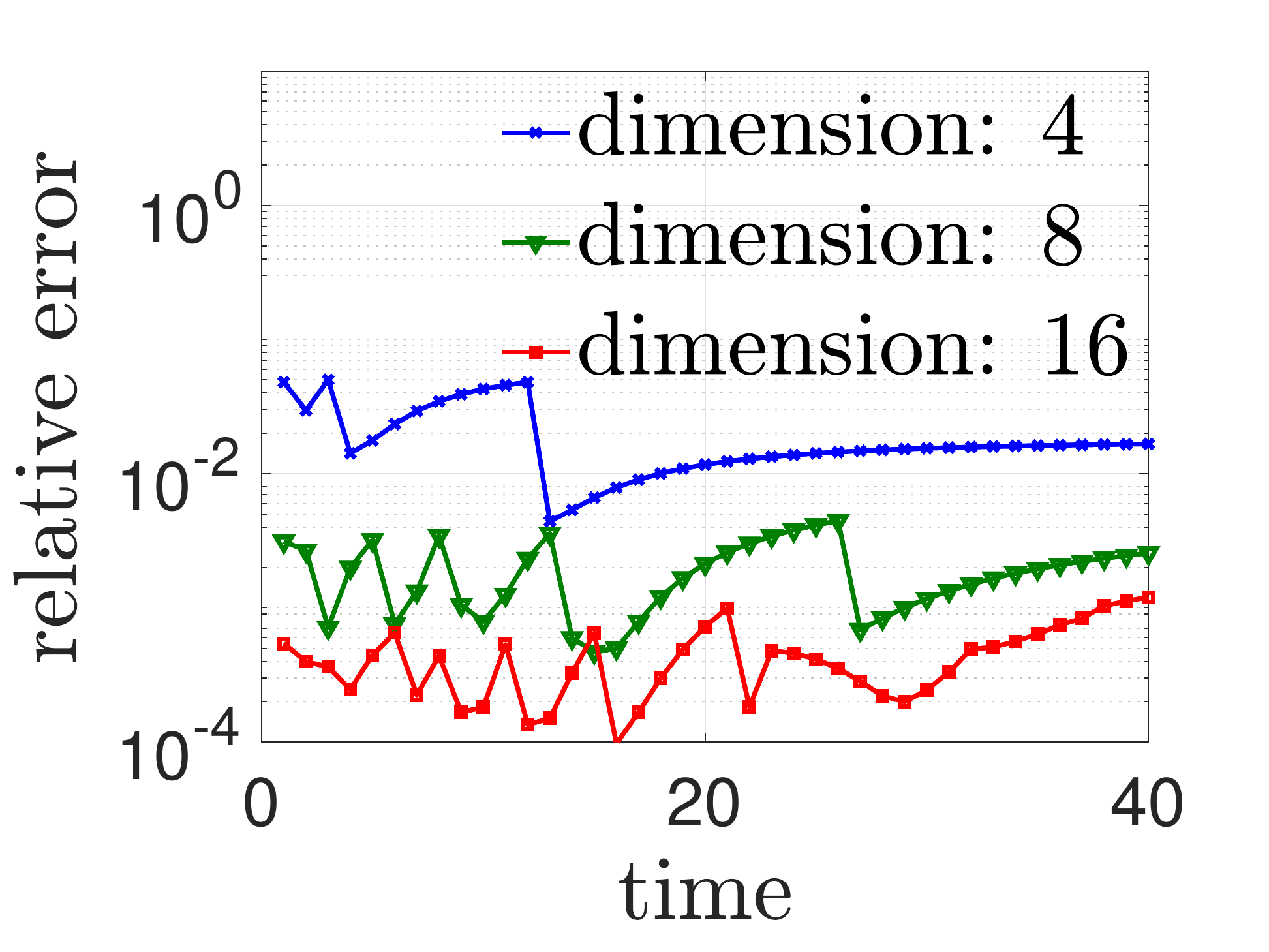}}
  \caption{Neutron flux and relative errors for the first problem} 
  \label{fig:example1}
  \end{figure}

  The space--time ROM is constructed, using $\svdtol = 2\times 10^{-8}$ and $\svtol =
  10^{-14}$ for the basis size of $16$, whose reduction factor is around
  twenty-seven million.  The ROM simulation uses $1$ core in RZTopaz.  With the
  basis size of $16$, the relative error with respect to the full order model
  solution is less than $0.1 \%$ as described in Figure~\ref{fig:example1}(c).
  Figure~\ref{fig:example1}(a) and (c) show the neutron flux distributions at
  the first and last time steps, respectively. The space--time ROM simulation
  with the basis size of $16$ takes $0.0055$ seconds, resulting in wall-clock
  time speed-up of $4,093.7$ and CPU time speed-up of $32,727.2$.

  \begin{figure}[h!]
  \centering
    \subfloat[Neutron flux at $t=0$ $nsec$\label{fig:initial2}]
    {\includegraphics[width=0.27\textwidth]{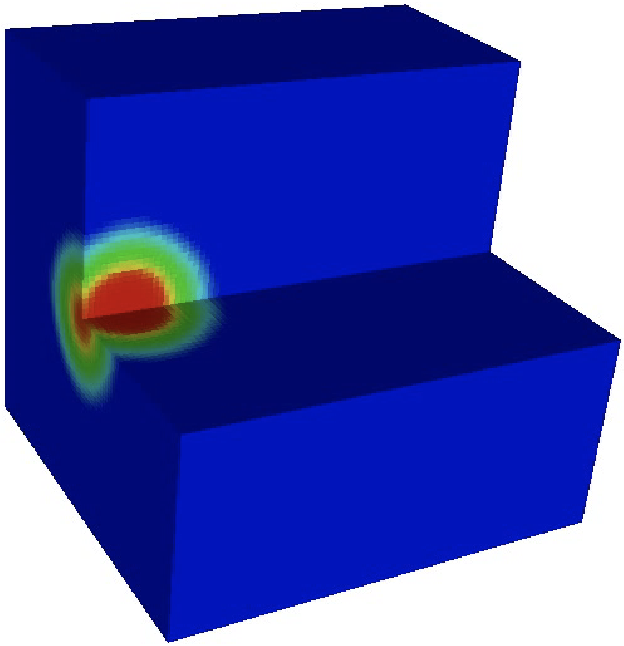}}\hfill
    \subfloat[Neutron flux at $t=30$ $nsec$\label{fig:final2}]
    {\includegraphics[width=0.27\textwidth]{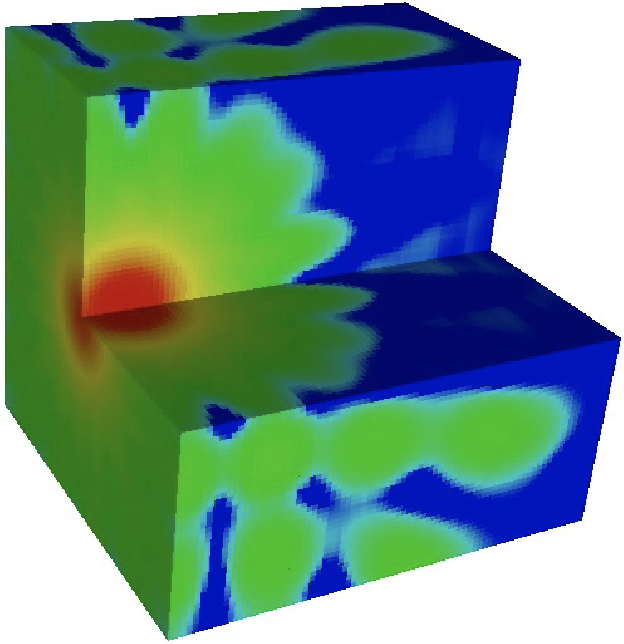}}\hfill
    \subfloat[relative error  \label{fig:relerror2}] 
    {\includegraphics[width=0.4\textwidth]{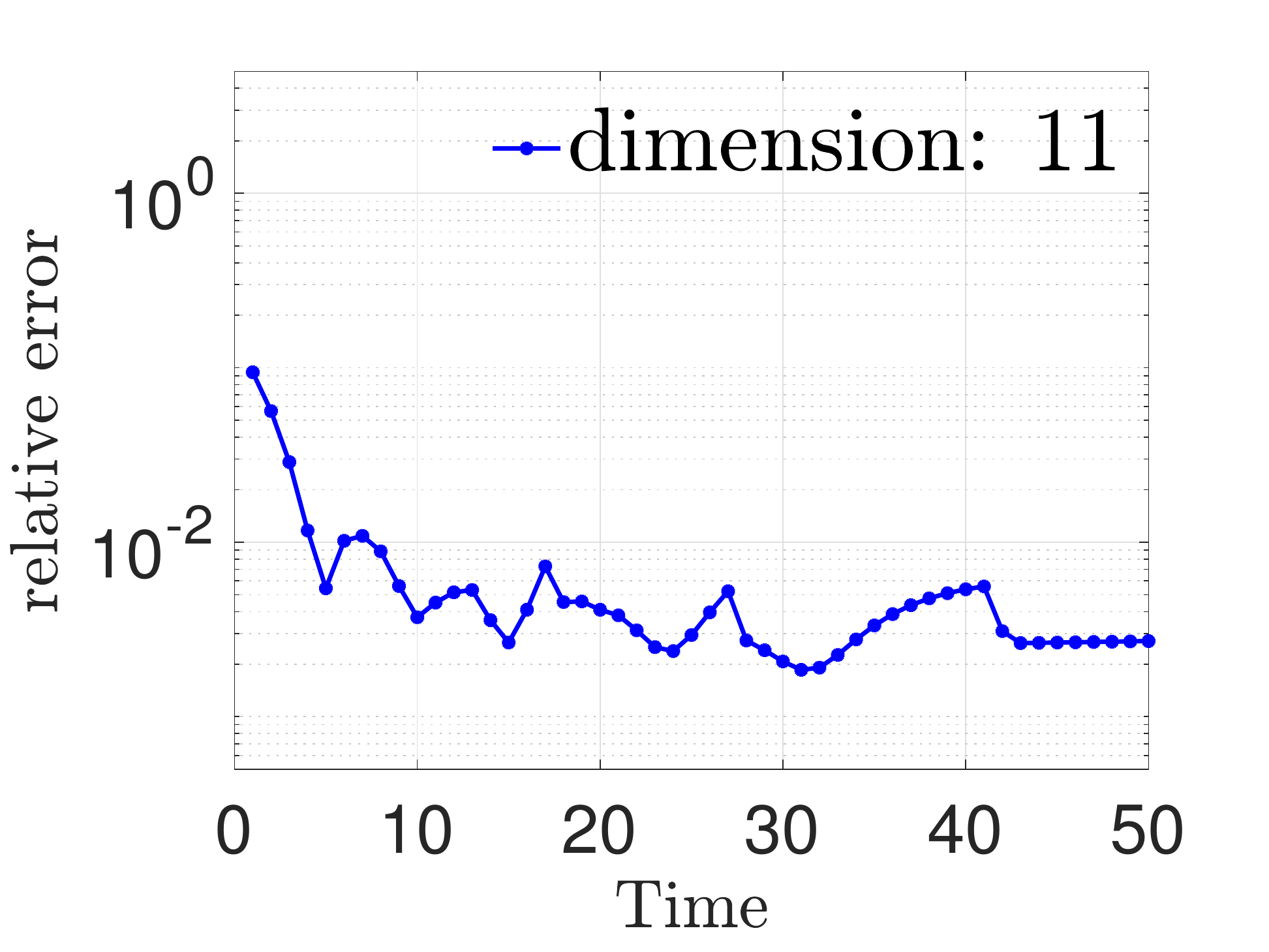}}
  \caption{Neutron flux and relative errors for the second problem, i.e., a truly 3D and larger-scale problem} 
  \label{fig:example2}
  \end{figure}

\subsection{Example 2: truly 3D case}\label{sec:ex2}
  The second example solves the Boltzmann particle transport equation on a
  different geometry, i.e., Figure~\ref{fig:ex2}. This is truly 3D with more
  structure than the previous example. The neutron source is $14.1$ $MeV$, which
  is in the 2nd energy group. The source is constant and the final simulation
  time step is at $\finaltime=30$ $nsec$ with a uniform time step $\timestep =
  0.2$ $nsec$.  The mesh is $80\times80\times80$, resulting in $512,000$ spatial
  zones.  There are eighty angular directions and seventeen energy groups.   As
  a result, there are $696,320,000$ degrees of freedom in space and
  $104,448,000,000$ degrees of freedom in space and time.  The full order model
  simulation uses $64$ cores in RZTopaz and takes 123.3 seconds, resulting in
  the CPU time of around 2.2 hours.

  The space--time ROM is constructed, using $\svdtol = 5\times 10^{-7}$ and $\svtol =
  10^{-14}$ for the basis size of $11$, whose reduction factor is around ten
  billion.
  The ROM simulation uses $1$ core in RZTopaz.  With the basis size of $11$, the
  relative error with respect to the full order model solution is less than $1.0
  \%$ as described in Figure~\ref{fig:example2}(c). Figure~\ref{fig:example2}(a)
  and (c) show the neutron flux distributions at the first and last time steps,
  respectively. The space--time ROM simulation with the basis size of $11$ takes
  $0.00682$ seconds, resulting in wall-clock time speed-up of $7,891.2$ and CPU
  time speed-up of $1,157,067.4$.

  \section{Conclusion}\label{sec:conclusion}
  Block structures in the space--time basis enable an efficient implementation
  of space--time reduced operators, which require small additional costs to the
  construction of the corresponding spatial ROM. Additionally, an incremental
  SVD is used to construct spatial and temporal bases in memory efficient way.
  As a result, the training cost of the space--time ROM is considerably reduced.
  Furthermore, because the space--time ROM achieves both space and time
  dimension reduction, considerably more reduction is accomplished than the
  spatial ROM, resulting in a great speed-up in online phase without losing much
  accuracy. It is demonstrated with Boltzmann transport problems where a
  reduction factor of twenty-seven million to ten billion and a CPU time
  speed-up of thirty-two thousand to one million were achieved by our
  space--time ROM.  Finally, our space--time ROM is {\it not} limited to a
  space--time full order model formulation. It is amenable to any time
  integrators although the backward Euler time integrator is used as an
  illustration purpose in this paper.

  Future works include applying the space--time ROM in the context of design
  optimization, uncertainty quantification, and inverse problems. Also, we will
  develop an efficient space--time ROM for nonlinear dynamical systems, such as
  TRT problems.

  \section*{Acknowledgement}\label{sec:acknowledgement}
  This work was performed at Lawrence Livermore National Laboratory and was
  supported by the LDRD program (17-ERD-026) and LEARN project (39931/520121).
  Lawrence Livermore National Laboratory is operated by Lawrence Livermore
  National Security, LLC, for the U.S. Department of Energy, National Nuclear
  Security Administration under Contract DE-AC52-07NA27344 and LLNL-JRNL-791966.

  \bibliographystyle{plain}
  \bibliography{references}

\end{document}